\magnification = \magstep1

\font\big=cmbx10 scaled \magstep2
\overfullrule=2pt
\def\part#1{\bigskip\bigskip\bigskip\noindent {\big #1}
            \bigskip\bigskip}

\def\section#1{\bigskip\medskip\goodbreak\noindent {\bf #1} \bigskip}

\def\kasten{\hfil\vrule height6pt width5pt depth-1pt\break}
\def\itemm{\smallskip\goodbreak
   \advance\nummer by1 \item{[\the\nummer] : }}
  \newcount\nummer
  \nummer = 0
\def\1{{1\!\!1}}
\def\3{\ss }
\def\AC{{\cal A}}

\def\BC{{\cal B}}
\def\b{\beta }
\def\C{\hbox{ \vrule width 0.6pt height 6.2pt depth 0pt \hskip -3.5pt}C}

\def\EC{{\cal E}}
\def\e{\varepsilon }
\def\FC{{\cal F}}
\def\f{\forall }
\def\G{\Gamma }
\def\GC{{\cal G}}
\def\g{\gamma }

\def\HC{{\cal H}}
\def\half{{1 \over 2}}
\def\IC{{\cal I}}

\def\L{\Lambda }

\def\l{\lambda }

\def\MC{{\cal M}}
\def\N{{I\!\! N}}

\def\PC{{\cal P}}

\def\R{{I\!\!R}}

\def\SC{{\cal S}}
\def\s{\sigma }

\def\TC{{\cal T}}

\def\Z{Z\!\!\! Z }

\def \b {\bigskip \noindent}
\def \m {\medskip \noindent}
\def \Ci{{C^\infty}}
\def \s {\smallskip \noindent}
 
\def \G {{\hbox{\bf  G} }}
\def \g {{ \bf \GC}}

\def \gld  {{\hbox{grad log det }}}

\def \H {{ \HC}}

\def \O {{ \Omega_e([0,1], G)}}
\def \HO {{ H^0([0,1], {\bf g})}}
\def \OO {{ \Omega_0([0, 1], {\bf g})}}

\centerline{\bf  REGULARISABLE AND MINIMAL
 ORBITS }
\m \centerline{\bf FOR GROUP ACTIONS IN INFINITE DIMENSIONS  }
\m  \b\centerline{\it M.Arnaudon and S.Paycha }\s
\centerline{\it Institut de Recherche Math\'ematique Avanc\'ee}
\centerline{\it Universit\'e Louis Pasteur et CNRS}
\centerline{\it 7, rue Ren\'e Descartes}
\centerline{\it 67084 Strasbourg Cedex France}\b\b
 \b {\bf Abstract}: We introduce a class of {\it regularisable}
 infinite dimensional  principal fibre bundles
  which includes fibre bundles arising in gauge field
theories like Yang-Mills and string theory and which generalise finite
dimensional Riemannian principal
fibre bundles induced by an isometric action. We show that
 the orbits of regularisable bundles
  have   well defined,  both heat-kernel and zeta function
regularised volumes. We introduce two notions of  minimality
 (which extend the finite dimensional
one) for these
orbits, using both heat-kernel and zeta function regularisation methods
and show they coincide. For each of these notions, we give an infinite
 dimensional version of Hsiang's theorem which extends the finite dimensional
case, interpreting minimal
 orbits as orbits with extremal (regularised) volume.
 \b {\bf 0. Introduction}\b
This article is concerned with the notions of  regularisability and minimality
of orbits for an isometric  action of an infinite dimensional Lie group $\G$ on
an infinite dimensional
 manifold $\PC$.  Our study is based on heat-kernel regularisation methods.
Notions of regularisability and minimality  have   already been studied by
other authors (see [KT], [MRT]) in a particular context and
using zeta function regularisation methods. We shall confront the different
approaches
to these notions as we go along.\s
We shall introduce a class of principal fibre bundles called
(resp. {\it pre-}){\it
 regularisable}
fibre bundles which generalise to the infinite dimensional case finite
dimensional Riemannian
principal fibre bundles arising from a free  isometric action.
We show that
the fibres of these  (resp. pre-)regularisable bundles
  have a well defined (both
heat-kernel and zeta function)  regularised (resp. preregularisable)
 volume which is G\^ateaux differentiable.
This class of (pre-) regularisable fibre bundles includes some infinite
dimensional principal bundles arising from gauge field theories
 such as Yang-Mills and string theory.  \s
We introduce various notions
of minimality, heat-kernel  minimality  and strong heat-kernel minimality using
  heat kernel regularisation methods on one hand and
zeta function minimality,
using  zeta function regularisation methods on the other hand, all of
which extend  the finite dimensional notion and coincide in the finite
dimensional
 case.
Whenever the structure group is equipped with a fixed Riemannian metric,
 we show that
 (strongly) minimal fibres of a (pre-)regularisable principal fibre bundle
coincide with  the ones
with extremal (pre-)regularised
 volume among orbits of the same type for the group action, the
regularisation being
taken in the heat-kernel sense.
This gives an infinite dimensional version of Hsiang's theorem on
(pre-) regularisable principal fibre bundles with structure group equipped
with a fixed Riemannian
 metric, which we extend (adding a term which reflects the
 variation of the metric on the structure group)
 to any  (pre- )regularisable principal bundle.
\m
  Starting from a systematic review of the notions of heat-kernel and
 zeta-function
regularised determinants in section I, in section II we introduce the notions
of
 regularisable principal
fibre bundle,
 heat-kernel ( pre) regularisability and  heat-kernel (strong)
minimality of orbits, relating (strong) minimality with
the G\^ateaux-differentiability of heat-kernel (pre-)regularised
determinants interpreted as volumes of
fibres. In section III,
we  compare these notions to zeta-function regularisability and
minimality, relating the  latter to
G\^ateaux-differentiability of zeta-function regularised determinants.
We show that     the two notions of regularisability
   and  minimality   coincide on the class of fibre bundles we consider.
  The relations we set up between the  regularised
mean  curvature vector and the directional gradients of the regularised
determinants
yield an infinite dimensional version of Hsiang's
theorem from both the heat-kernel  and the zeta function point of view.
In  Appendix A , we   apply these   results
to
 the coadjoint action of a loop group  thus recovering some results
concerning regularisability and minimality of fibres studied in
  [KT]. In Appendix B, we investigate     minimality of the orbits in the
case of Yang-Mills
action, for which  a notion of (zeta function )  minimality had been
suggested in  [MRT] from which our
 notion differs slighlty.
 We point out the fact that when  the underlying manifold is
of dimension 4, only  if the irreducible
connections are Yang-Mills,
 do  the notion of    zeta function   and heat-kernel  minimality coincide.
    In both   examples, the space $\PC$, resp.
 the group $\G$ are modelled on a space of sections
of a vector bundle $\EC$, resp. $\FC$ with finite dimensional fibres  on a
closed finite dimensional
manifold $M$ and $\G$ acts on $\PC$ by isometries. \s
One could show, in  a similar way to the Yang-Mills case, using results of [RS]
that the bundle $\MC_{-1}\to \MC_{-1}/\hbox{Diff}_0$ (described in [FT])
arising in bosonic string theory ( where $\MC_{-1}$ is the manifold of smooth
Riemannian metrics with curvature $-1$ on a compact boundaryless
Riemannian surface of genus greater than 1
  and $\hbox{Diff}_0$ is the group of smooth diffeomorphisms of the
 surface which are homotopic to zero), is also a regularisable fibre bundle
so that most results of this paper can be applied to   this fibre bundle.
  However since, unlike the case of Yang-Mills theory, its structure group
  $ \hbox{Diff}_0$ is not equipped with a fixed Riemannian structure but
    with a  family of Riemannian metrics which is parametrised by
 $g\in\MC_{-1}$, minimality of the fibres is not equivalent to
extremality of the volumes of the fibres (see Proposition 2.2),
and we chose not to treat this example in detail in this paper. \m
 The geometric notions developped in this paper play a important role
 when projecting  a class of semi-martingales
defined on the total manifold onto the orbit space for a certain class of
 infinite dimensional group actions. The heat-kernel regularisation method
 yields natural links between
 the geometric and the stochastic picture, which we investigate in [AP2].
 The stochastic picture described in [AP2] leads to a stochastic
interpretation of the Faddeev-Popov
 procedure used in gauge field theory to
 reduce a formal volume measure on path space to a measure on the
orbit space, the formal density of which is a regularised "Faddeev-Popov"
determinant.
\b{\bf Acknowledgments}: We would like to thank Steve Rosenberg most warmly
for very valuable
 critical
 comments  he made on a previous version of this paper. We would also like
to thank David Elworthy for his generous hospitality at the Mathematics
 Department
of Warwick University where part of this paper was completed.  \b\b
\b { \bf I. Heat-kernel and Zeta-function regularized determinants}
\b
 In this section, we recall some  basic facts about heat-kernel
 and zeta function regularised determinants, comparing the two regularisations.
Although the results presented here are well known and  frequently used in
the physics
 literature,
  it
 seemed necessary to us to give a
clear and precise presentation of  the heat-kernel and
 zeta function regularisation  procedures for later use.
 \s
Let us first introduce some notations.
 For a function $t\mapsto f(t)$,
 defined on an interval  of $\R^{+*}$ containing $]0,1]$,
we shall write
$f(t)\simeq_0 \sum_{j=-J}^{K-1} a_j t^{j\over m}$, $a_j\in \R$,
$J, K\in \N$, $m\in \N^*$,  if
there exists a constant $ C >0$, such that
 $$ \vert f(t)- \sum_{j=-J}^{  K-1} a_j t^{j\over m}\vert <C  t^{K \over m}
 \quad \f
\quad 0<t<1\eqno(1.0)$$  \s
In the following, we shall always assume that $J\geq m$.\b
{\bf Lemma 1.0}:
Let $(A_\e),\e\in ]0,1] $ be a one parameter family of trace-class
operators on a separable
Hilbert space $H$ (in particular tr($A_1$) is finite)
such that:
\item{1)}$\e \to \hbox{tr} A_\e$ is differentiable on $]0,1[$
\item{2)} $\exists J   \in \N$,$m\in \N^*$, $(a_j)_{j\in \{-J, \cdots,
-1\}}, a_j\in \R$ such that
$${d\over d\e} (\hbox{tr}A_\e)\simeq_0\sum_{j=-J}^{-1}\e^{j\over m} a_j,
      \eqno (1.1)$$
 Then     the expression
 $\displaystyle \hbox{tr}A_\e-\sum_{j={-J+m}   }^{ -1}
{ma_{j-m}   \over j}\e^{j\over m}- a_{-m} \hbox{log} \e$ converges
 when $\e \to 0$.
 \m {\bf Remark }: In the following, we shall not distinguish the two cases
and adopt the convention that
the sum from $-J+m$ to $-1$ is zero whenever $J=m$.
\m {\bf Proof}:
To show this, let us set for $0<\e<1$, $g_\e= \hbox{tr}A_\e-
\sum_{j={-J+m}, j\neq 0}^{ -1+m}
{ma_{j-m}   \over j}\e^{j\over m}-a_{-m} \hbox{log} \e$.
  For
 $0<\e<\e^\prime<1$, we have:
 $$\eqalign{ \vert  g_\e -g_{\e^\prime}\vert
  & \leq \int_\e^{\e^\prime} \vert
 {d\over dt} (\hbox{tr}A_t)-\sum_{j=-J}^{ -1}t^{j\over m} a_j\vert dt\cr &\leq
{C  }  ( \e^\prime-\e )  \leq C  \e^\prime\quad \hbox{ by (1.1)}\cr}$$
so that $(g_\e)$ is a Cauchy sequence and hence converges when $\e \to 0$.
  From the convergence of $g_\e$ then follows the
 convergence of $ \hbox{tr}A_\e-\sum_{j={-J+m},j\neq 0}^{ -1}
{ma_{j-m}   \over j}\e^{j\over m} -a_{-m}\hbox{log} \e$  when $\e$ goes
to zero since the terms indexed by $0,1, \cdots, -1+m$ converge when $\e \to
0$.

  \kasten
\b {\bf Definition}: Whenever   $ \hbox{tr}A_\e-\sum_{j={-J+m}  }^{ -1}
{ma_{j-m}   \over j}\e^{j\over m} -a_{-m} \hbox{log} \e$ converges
when $\e \to 0$. We
  shall
 call the limit  the {\it  regularized limit trace} of the family $\AC=
(A_\e)  $
 and denote
 it by
 $\hbox{tr}_{reg} (\AC)$ so that
$$\hbox{tr}_{reg} (\AC)=\lim_{\e \to 0} (\hbox{tr}A_\e-\sum_{j={-J+m}}^{ -1}
{ma_{j-m}   \over j}\e^{j\over m}-a_{-m} \hbox{log} \e )
 \eqno (1.2)$$
This regularised limit trace depends of course on   the whole  one parameter
family and
  on the choice of the parameter $\e$.
 \b\b Let $\BC=({B_\e}) $  be a one parameter family of strictly
 positive self adjoint
 operators such that
$\AC= (A_\e)$ with $A_\e\equiv\hbox{log} B_\e$ is a family of trace class
operators. We can define  the determinant of $B_\e$ as $\hbox{det} B_\e=
e^{ \hbox{tr log} B_\e}$.
If the family $\AC=(A_\e)$ has a
 regularized limit trace,  we shall define the {\it
regularized limit} determinant of the family $\BC$ by
$$\hbox{det}_{reg} \BC\equiv e^{\hbox{tr}_{reg} (\AC)}\eqno(1.3)$$\b
We now introduce a family of heat-kernel operators which play a
fundamental role in this paper. For this we define  for $\e>0$ a
function  $ h_\e$  by:
$$\eqalign{h_\e:\R^{+*}&\to \R\cr
 \l&\mapsto e^{-\int_\e^\infty {e^{-t\l}\over t }dt}\cr}$$
Notice that $h_\e$ is $ C^\infty$, non decreasing  and
$ (\hbox{log}h_\e)^\prime(\l)=\l^{-1}e^{-\e \l}$.
Writing $\hbox{log}h_\e(\l) -\hbox{log} \e=
-\int_{\e  }^\infty {e^{-t}\over t}
 dt -\int_\e^{\e \over \l}
{ e^{-\l t }\over t} dt+\int_\e^1 {1\over t}dt$, we find
that $$\lim_{\e \to 0}{ h_\e (\l)\over \e} =\l e^{\int_0^1
{1-e^{-t} \over t} dt -\int_1^\infty  {e^{-t} \over t}
 dt} \eqno (1.3 bis)$$
For a strictly positive self adjoint operator $B$ on a Hilbert space $H$,
we can define
$h_\e(B)$ which yields a one parameter family of operators
$B_\e\equiv h_\e(B)$.
\m {\bf Definition}: If $\hbox{log}h_\e(B)$ is trace class, we can define
$$\hbox{det}_\e(B)= e^{\hbox{tr log} h_\e(B)}\eqno(1.4)$$
  \m
 {\bf Definition}:  Let $B$ be a strictly positive self-adjoint operator
on a separable Hilbert space. Whenever the one
parameter family $\BC=(h_\e(B))$
has a regularized limit determinant, we shall call this limit the {\it
heat-kernel
regularized determinant} of $B$ and we denote it by
$\hbox{det}_{reg}(B)$.\s
We have
$$\hbox{det}_{reg}B= \hbox{det}_{reg} (\BC)= e^{\hbox{tr}_{reg} (\AC)}
\eqno (1.5)$$
with
$\BC=(h_\e(B))$, $\AC=(\hbox{log} h_\e(B))$.\s
 In the
following, we give conditions under which we can define the heat-kernel
regularized
determinant of an operator $B$. But before that, let us state an easy lemma
which will prove to be useful for what follows.\m
\b
{\bf Lemma 1.1}:Let $B$ be a strictly positive self-adjoint operator
on a separable Hilbert space
 such that
  $e^{-\e B}$  is trace class for some $\e>0$. Then
$A_\e\equiv \hbox{log} h_\e(B)$ is also trace-class.
\m {\bf Proof}: Since $e^{-\e B}$ is trace class, it is compact and hence
has a purely discrete spectrum $\{\mu_n, n\in \N\}$, $\mu_n>0$, $\mu_n$
tending to zero. Hence  $ B=-\e^{-1}\hbox{log }e^{-\e B}$ also has purely
discrete spectrum
$\l_n=-\e^{-1} \hbox{log}\mu_n$ which tends to infinity and is bounded from
below by a
strictly positive constant.
Let $\L_{n_0} $ be a non zero eigenvalue and let us index the eigenvalued in
encreasing order so that   $ \l_n \leq \l_{n +1} $.
We have:
$$\eqalign{ \sum_{n\geq n_0} \int_\e^\infty {e^{-t\l_n} \over t}dt &\leq
\e^{-1} \sum_{n \geq
n_0}\l_n^{-1} e^{-\e \l_n}\cr
&\leq \e^{-1} \l_0^{-1} \sum_{n\geq n_0} e^{-\e \l_n} \cr}$$
The last expression is finite by assumption so that
$  \sum_{n\geq n_0}\vert\int_\e^\infty {e^{-t\l_n}\over t} dt\vert $ is
finite and $A_\e$ is trace class.\kasten
\b
\m {\bf Lemma 1.2 }: Let $B$ be a strictly positive self-adjoint operator
on a separable Hilbert space
 such that
\item{1)} $e^{-\e B}$ is trace class for any $ \e>0$.
 \item{2)} There is a family $(b_j)_{j=-J,\cdots, 0}$ , $b_j\in \R$ and an
integer $m>0$
such that
$$\hbox{tr} e^{-\e B}\simeq_0\sum_{j=-J}^{m-1} b_j \e^{j\over m}.  $$
    \s Then the operator $B$ has a   heat-kernel regularised determinant
and   we have
$$\eqalign{\hbox{det}_{reg}B&= \lim_{\e \to 0}\left(\hbox{det}_\e B
 e^{ -\sum_{j=-J}^{-1} {mb_j\over j}\e^{j\over m} -b_0 \hbox{log} \e   }\right)
\cr
&= e^{ \left(-\sum_{j=-J,j\neq 0 }^{m-1}
 {mb_j\over j}-\int_1^\infty \hbox{tr}
{ e^{-tB} \over t}dt -\int_0^1 {F (t)\over t} dt\right)}\cr}\eqno (1.6)$$
with
$$F(t)= \hbox{tr} e^{-tB}-\sum_{j=-J}^{m-1} b_j t^{j\over m}.\eqno (1.7)$$
  \m {\bf Remark}: If the Hilbert  space $H$ the operator $B$
acts on is finite dimensional
of dimension $d$, since
$\lim_{\e \to 0}\hbox{tr}e^{-\e B}=d= b_0$,   (1.6) yields
 $  \hbox{det}_{reg} B=\lim_{\e \to 0}( \hbox{det}_\e B\e^{-d})$.
  \m {\bf Proof}:
By Lemma 1.1, we know that
$A_\e=-\int_\e^{\infty} {e^{-tB} \over t} dt= \hbox{log} h_\e(B)$
 is trace-class.
Since all the terms involved are positive, we can exchange the integral and sum
symbols so that
  $\hbox{tr} A_\e= -\int_\e^\infty \hbox{tr} {e^{-tB} \over t} dt$.
We now apply lemma 1.0 to show that the family $A_\e$ has a regularized
limit trace.
The map
$t  \to \hbox{tr} A_t$ is differentiable and by assumption 2), for $t>0$,
we have
 $$ {d\over d  t} \hbox{tr} A_t  = \hbox{tr} {e^{-t
B} \over t} = \sum_{j=-J}^{ m -1 }
  b_j t^{j-m\over m} +{F (t)\over t}= \sum_{j=-J-m}^{-1} b_{j+m} t^{j\over
m} +{F(t) \over t}\eqno (1.8)$$
where $F$ is as in (1.7) and $\vert F(t)\vert \leq Ct$.
 Thus, taking $K=0$ in (1.0), we have:
$$ {d\over d  t} \hbox{tr} A_t  \simeq_0   \sum_{j=-J-m}^{ -1 }
  b_{j +m}t^{j \over m} $$
so that we can apply
  Lemma 1.0 from which follows (replacing $J$ by $J+m$ and $a_j$ by
$b_{j+m}$) )
that  the one parameter family  $\AC=(A_\e)$ has a
finite regularized
 limit trace
$\hbox{tr}_{reg}(\AC)\equiv \lim_{\e\to 0}
(\hbox{tr} A_\e-\sum_{j=-J+m}^{-1} {mb_j\over j}\e^{j\over m}
-b_0 \hbox{log} \e)
 $.   Hence, by (1.5)   $B$  has  a heat-kernel   regularized   determinant:
$$\eqalign{\hbox{det}_{reg}B=e^{\hbox{ tr}_{reg} \AC}&=\lim_{\e \to 0}
\left( e^{\hbox{tr} A_\e}
 e^{-\sum_{j=-J}^{-1} {mb_j\over j} \e^{j\over m} -b_0
\hbox{log}\e }\right)\cr
&=
 \lim_{\e \to 0}
\left( \hbox{det}_\e B
 e^{-\sum_{j=-J}^{-1} {mb_j\over j} \e^{j\over m} -b_0
\hbox{log}\e }\right).\cr}$$
Since
$\hbox{tr} A_1=-\int_1^\infty \hbox{tr} {e^{-tB} \over t} dt$,
 integrating (1.8) between $\e$ and $1$ yields
$$ \hbox{tr} A_\e-\sum_{j=-J, j\neq 0}^{m-1} {mb_j\over j}
 \e^{j\over m} -b_0 \hbox{log} \e =  -\sum_{j=-J, j\neq 0}^{ m -1} {mb_j
\over j}
-\int_\e^1 {F (t) \over t} dt-\int_1^\infty \hbox{tr} {e^{-tB} \over t}dt.
\eqno (1.9)  $$
Since $m\geq 1$, we have
$$\lim_{\e \to 0}( \hbox{tr} A_\e-\sum_{j=-J}^{ -1} {mb_j\over j}
 \e^{j\over m} -b_0 \hbox{log} \e)= \lim_{\e \to 0}( \hbox{tr} A_\e-
\sum_{j=-J, j\neq 0}^{m -1} {mb_j\over j}
 \e^{j\over m} -b_0 \hbox{log} \e)$$ which
  combined with (1.9)  and using (1.4) yields  (1.6). \s
Notice that the integral
$\int_\e^1 {F (t) \over t} dt $ is absolutely convergent.
 Indeed, by assumption 2), $\vert F(t)\vert \leq C t$ for some positive
constant $C$.
 \kasten
   \b
The following lemma gives a class of operators which
 fit in the framework described above.\b
{\bf Lemma 1.3}:
Let $B$ be a strictly positive self adjoint
 elliptic operator of order $m>0$ on a compact boundaryless manifold. For
any $\e>0$,
  $e^{-\e B}$ is  trace class and $B$  has a well
defined
 heat-kernel regularised   determinant. \m {\bf Proof}:
 We shall   show that the assumptions of
Lemma 1.2 are fulfilled.
\s Condition  1) in Lemma 1.2 follows from the fact that
a strictly positive s.a  elliptic operator on a compact boundaryless
manifold has purely discrete spectrum
$(\l_n )_{n\in \N}$, $\l_n>0$, $\l_n \simeq Cn^\alpha$, for some $C>0$,
 $\alpha >0$ ( see e.g [G] Lemma 1.6.3).
Indeed, from this fact easily follows that $\hbox{tr}e^{-\e B}=
\sum_n e^{-\e \l_n}$ is finite.
\s Conditions 2) of Lemma 1.1 follow from the fact that for
  a s.a elliptic operator  $B$ of order $m$ on a compact manifold  of
dimension d without
boundary,
$\hbox{tr} e^{-tB}\simeq_0 \sum_{ j=-d}^{K-1} a_jt^{j\over m} $      for
any $K>0$
(this follows for example from
 Lemma 1.7.4 in [G]).
  Applying lemma 1.0, we can therefore define the heat-kernel regularized
determinant of $B$.
 \kasten\m
 The above definition extends to a class of positive self-adjoint operators
which satisfy requirements
 1) and 2) of Lemma 1.2 and have possibly non zero kernel. Requirement  1)
of the
 lemma implies that this kernel is finite dimensional. Let $P_B$ the
orthogonal projection
 onto the kernel of the operator $B$ acting on $H$ and let us set
$H^\perp\equiv
 (I-P_B)H$. Let us assume that $H^\perp$ is invariant under the action of
 $B$ so that we can consider the  restriction  $B^\prime\equiv B/H^\perp$.
 The operator $B^\prime$ satisfies requirements of Lemma 1.2, namely
 \item{1)} $e^{-\e B^\prime}$ is trace class for any $ \e>0$.
 \item{2)} There is a family $(b_j^\prime)_{j=-J,\cdots, 0}$ ,
$b_j^\prime\in \R$ and an integer $m>0$
such that
$$\hbox{tr} e^{-\e B^\prime}\simeq_0\sum_{j=-J}^{m-1} b_j^\prime \e^{j\over
m}.  $$
where $b^\prime_j=b_j$ for $j\neq 0$ and
$b_0^\prime=b_0-dim Ker B$.
    \s Under
     assumption 1), we can extend definition (1.4) and define:
     $$\hbox{det}_\e^\prime B\equiv e^{\hbox{tr log} h_\e(B^\prime)}$$
     Under   assumptions 1) and 2), the operator $B^\prime$ has a
       heat-kernel regularised determinant
 $$\eqalign{\hbox{det}_{reg}^\prime B&\equiv \lim_{\e \to 0}\left(\hbox{det}_\e
 B^\prime
 e^{ -\sum_{j=-J}^{-1} {mb_j\over j}\e^{j\over m} -b_0 \hbox{log} \e
 +(\hbox{ dim Ker B})\hbox{log} \e}\right)
\cr  } $$\b\b
 Let us at this stage compare the heat-kernel regularised determinant with the
zeta-function regularised one. We refer the reader to [AJPS], [G]
 for a precise description of the zeta-function regularisation procedure
and only describe the main lines of this procedure here.
 \s
 Recall that for a strictly positive self adjoint operator $B$ acting on
 a separable Hilbert space with purely discrete spectrum given by the
eigenvalues $(\l_n, n\in \N)$ with the property $   \l_n\geq C n^\alpha, C>0,
 \alpha >0$ for large enough $n$, we can define the zeta function of $B$
 by:
$$\zeta_B(s)\equiv \sum_n \l_n^{-s}, \quad s\in \C, \quad \hbox{Re}s>
{1\over \alpha}$$
Furthermore, $\zeta_B(s)$ admits a meromorphic continuation on
the whole plane  (see e.g [G] Lemma 1.10.1)
which is regular at $s=0$ and one can define the zeta function regularized
determinant of $A$ by
$$ \hbox{Det}_{reg}(B)= e^{-\zeta_B^\prime(0)}\eqno (1.10)$$
\m {\bf Remark}: From the definition, easily follows that
in the finite dimensional case the zeta-function regularised and the
ordinary determinants coincide.
\m
 \m The following lemma compares the  two regularizations.
\b {\bf Lemma 1.4}: Let
$B$ be  a strictly positive self-adjoint densely
 defined operator on  a Hilbert space $H$ such that
 \item{1)}$B$ has purely discrete spectrum $(\l_n)_{n\in \N}$ with
 $\l_n\geq Cn^\alpha$,$ C>0, \alpha>0$ for large enough $n$,
 \item{2)}
$$\exists J>0, m>0, (b_j)_{j=-J,\cdots,m-1} \quad \hbox{ such that}\quad
\hbox{tr} e^{-\e B}\simeq_0 \sum_{j=-J }^{m-1} b_j \e^{j\over m}$$
\item{} Then
$$\eqalign{\hbox{Det}_{reg}B&= e^{ -\gamma b_0} ( \hbox{det}_{reg} B)
\cr
&=e^{-(\gamma b_0+\sum_{j=-J,j\neq 0}^{m-1}{mb_j\over j}+\int_1^\infty {tr
e^{-tB }\over t}dt +\int_0^1 {F(t)\over t}dt)} \cr} \eqno(1.10 bis)$$
 where   $\gamma =\lim_{n\to \infty} (1+\half+\cdots +
{1\over n}-\hbox{log}n)$ is the Euler constant
 and $b_0$ is the coefficient arising in the heat-kernel expansion of $B$,
$$\hbox{tr}e^{-\e B}\simeq_0 \sum_{j=-J}^{m-1} b_j \e^{j\over m}$$ for some
 $J\in \N$, $m>0$.
\m {\bf Remark }:
   A proof of this result for
 the Laplace operator on a compact
Riemannian surface without boundary
can be found in [AJPS].
  \m {\bf Proof}:
Before starting the proof, let us recall that the function Gamma is defined by
$\Gamma(z)=\int_0^\infty {e^{-t}\over t} t^z dt$ for $0<\hbox{Re}z$.
Moreover $\Gamma(z)^{-1}$ is an entire  function and we have
$$\Gamma(z)^{-1}=z e^{\gamma z} \prod_{n=1}^\infty (1+{z\over n})
e^{-z\over n}$$
where $\gamma$ is the Euler constant. From this follows that in a
neighborhood of zero, we have the asymptotic expansion
$\Gamma(s)^{-1}=s +\gamma s^2 +O(s^3)$.\m
 Using the Mellin transform of the function $$\l^{-s}=
\Gamma(s)^{-1} \int_0^{+\infty} t^{s-1} e^{-t\l} dt$$
 we can write:
$$\Gamma(s) \zeta_B(s)=\int_0^1 t^{s-1} \hbox{tr} e^{-tB} dt +
 \int_1^\infty   t^{s-1} \hbox{tr} e^{-tB}dt\eqno (1.11)$$
Notice that the last expression on the r.h.s converges for $\hbox{Re}s\leq R$,
$R>0$
for, setting  $C_R=\hbox{sup}_n \hbox{sup}_{t \geq 1} t^{R-1}e^{-\half
t\l_n}$, we have
$\int_1^\infty t^{R-1} e^{-t\l_n}\leq C_R
\int_1^\infty e^{-\half t\l_n}=2C_R \l_n^{-1} e^{-\half \l_n}$
which is the general term of a convergent series.\s
As before we set
$$F (t)\equiv \hbox{tr} e^{-tB}-\sum_{j=-J}^{m-1} b_jt^{j\over m} \eqno
(1.12)$$
  Using (1.11 ) and (1.12), we can write for $s\in \C$
with large enough
 real part, $\hbox{Re}s>{J\over m}$:
$$\zeta_B(s)=  \Gamma(s)^{-1} \left(\sum_{j=-J}^{m-1} {b_j \over {j\over m}+s}+
 \int_1^\infty   t^{s-1} \hbox{tr} e^{-tB}dt+ \int_0^1 t^{s-1}F (t)dt
\right)
$$
  This equality then extends to an equality of meromorphic functions on
$\hbox{Re}s>0$ with poles $s={-j\over m}$.
  Using the
asymptotic expansion of
  the inverse of the Gamma function $\Gamma(s)^{-1}$  around zero,
   we have:
    $$\eqalign{ &\zeta_B^\prime(s) =
( 1+2\gamma s +O(s^2)) \left(
\sum_{j=-J}^{m-1 } {b_j \over {j\over m}+s} +
\int_1^{\infty}  t^{s-1} \hbox{tr} e^{-tB}dt+ \int_0^1t^{s-1} F (t)
dt\right)\cr
&+ (s+\gamma s^2 +O(s^3))
 \left(-\sum_{j=-J}^{m-1} {b_j \over ({j\over m}+s)^2} +\int_0^1  t^{s-1}
F (t)
\hbox{ln}
 (t) dt+ \int_1^\infty   \hbox{ln} (t) t^{s-1} \hbox{tr} e^{-tB}dt
\right)\cr}$$
Letting $s$ tend to zero, $s>0$,  since the divergent terms ${b_0\over s}$ and
$-s {b_0\over s^2}$
arising in each of the terms of this last sum compensate, we get:
 $$\zeta_B^\prime(0)=
  b_0 \gamma +  \left( \sum_{j=-J,j\neq 0 }^{ m-1} {mb_j \over j}
+ \int_0^1 {F (t)\over t} dt +\int_1^\infty
 { \hbox{tr} e^{-tB}\over t}dt\right)   $$
Hence, comparing with the expression of $\hbox{det}_{reg}B$ given in (1.6), we
find:
$$\hbox{log  Det}_{reg} (B)=-\zeta_B^\prime(0)=-b_0 \gamma +  \hbox{log det
}_{reg}(B)$$
and hence the equality of the lemma. \kasten
\b {\bf Remarks}: \item{1)} Notice that the same proof as in the lemma
replacing the function $\zeta_B(s)$ by
$\zeta_\l(s)=\l^{-s}$ for $\l>0$ (which boils down to taking a one
 dimensional space $H$ and $F(t)= e^{-t\l}-1$)
yields
$$-\hbox{log} \l= \gamma+
\int_0^1{e^{-t\l}-1\over t}dt +\int_1^\infty dt {e^{-t\l} \over t} $$
and when choosing $\l=1$, the integral representation of the Euler constant.
$$\gamma=
 \int_0^1{1-e^{-t } \over t}dt -\int_1^\infty dt {e^{-t } \over t}.$$
\item{2)} In the finite dimensional case, dim$H=d$, since $\lim_{\e \to 0}
\hbox{tr}e^{-\e B}=d=b_0$, from  the result of lemma 1.4  and the fact that
the zeta function regularised determinant coincides with the ordinary one,
follows that
 $$\hbox{det}_{reg}B=e^{ d\gamma }\hbox{Det}_{reg} B=
 e^{  d\gamma}\hbox{det} B=\lim_{\e \to 0} (\hbox{det}_\e(B) \e^{-d}) \eqno
(1.13)$$
 where detB denotes the ordinary determinant of $B$. This agrees with
(1.3 bis) using the integral representation given above of the Euler constant.
  \item{3)}Let  $M$ be a Riemannian manifold
     of dimension $d$
 and $B$ a positive self-adjoint elliptic
operator  with smooth coefficients acting on sections of a vector bundle
$V$  on $M$ with finite dimensional
 fibres of dimension $k$.  We know by [G] Theorem 1.7.6 (a)  that
$b_0= 0$ if $n$ is odd. However, in general  the coefficient $b_0$ is a
 complicated expression  given  in terms of the jets of the symbol of the
operator B.
In the following we shall be concerned with
the dependence of $b_0$ on the geometric data
given on that manifold.
\b   The notion of $\zeta$ function regularized determinant therefore
extends to an operator satisfying assumptions of Lemma 1.2, by setting:
$$\hbox{Det}_{reg} B=
 e^{   \gamma b_0}\hbox{det}_{reg} B \eqno(1.14)$$
 It furthermore extends to positive operators satisfying assumptions of Lemma
 1.2 with non zero (finite dimensional) kernel and which leave its
orthogonal supplement
 invariant, for in that case, we can set:
$$\hbox{Det}_{reg}^\prime B\equiv
 e^{   \gamma( b_o-dim Ker B)}\hbox{det}_{reg}^\prime B \eqno(1.15)$$
\b\b {\bf II  Regularisable  principal fibre bundles}\s
The aim of this section is to define a class of principal fibre bundles
for which we can define a notion of  regularised volume of the fibres and
for which
 these regularised volumes have
   differentiability properties.
\b  Let $\PC$ be a Hilbert   manifold equipped with a
(possibly weak) right invariant  Riemannian structure. The scalar product
induced on
 $T_p\PC$ by this Riemannian structure will be denoted by $<\cdot, \cdot>_p$.
 We shall assume this Riemannian structure induces
a Riemannian connection denoted by $\nabla$ and an exponential map with
the usual properties. In particular, for all $p_0$, $\hbox{exp}_{p_0} $
yields a
   diffeomorphism of a neighborhood of $0$ in
  the tangent space $T_{p_0}\PC$ onto a neighborhood of $p_0$ in
 the manifold $\PC$.
  \b Let   $\G$ be a     Hilbert
   Lie group ( in fact a right semi-Hilbert Lie group i.e a Hilbert Lie group
 in the usual sense up to the fact that only right multiplication is
required to be smooth in the sense of [P] is enough here)
acting smoothly on $\PC$ on the right
 by an isometric action
$$\eqalign{\Theta: \G \times \PC &\to \PC \cr
(g,p)&\to p\cdot g\cr}\eqno(2.0)$$
Let for $p\in \PC$
$$\eqalign{\tau_p:\GC&\to T_p\PC\cr
u&\mapsto {{d\over dt } (p\cdot e^{tu})}_{t=0}\cr}\eqno(2.0 bis)$$
where  $\g$ denotes the Lie algebra of $\G$.\s
We shall   assume that
  the action $\Theta$ is free (so that $\tau_p$ is injective on $\GC$) and that
  it induces a smooth manifold structure on
 the quotient space $\PC/\G $ and a smooth principal fibre bundle structure
given by the canonical projection $\pi:\PC\to \PC/\G $.
  \m
 Let us  furthermore equip the group $\G$  with a smooth family of
equivalent (possibly
weak)  $\hbox{Ad}_g$ invariant Riemannian metrics indexed by
$ p\in \PC$.
 The scalar product induced on $\g$ by the Riemannian metric on $\G$
indexed by $p\in \PC$ will be denoted by $(\cdot, \cdot)_p$. Since
the metrics are all equivalent, the closure of
$\g$ w.r.t $(\cdot, \cdot)_p$ does not depend on $p$ and we shall denote
 it by $H$.
  \m Since $\g$ is dense in $H$, $\tau_p$ is a densily defined
 operator on
   $H$ and
 we can define its adjoint operator $\tau_p^*$ w.r. to
the scalar products $(\cdot, \cdot)_p$ and $<\cdot, \cdot>_p$.
\s We shall
 assume  that $\tau_p^*\tau_p$ has a
 self adjoint extension on a dense domain $D(\tau_p^* \tau_p)$ of $H$.
\m We shall assume that $\tau_p$ is injective.
\s {\bf Warning}: Although $\tau_p$ is injective on $\GC$, the operator
$\tau_p^*\tau_p$
 might not be injective
on   the domain $D(\tau_p^* \tau_p)$ as we shall see in applications
(cfr.Appendix
A).
 \b {\bf Definition}: The orbit of a point $p_0$ is {\it   volume
preregularisable}
if   the following   assumptions
 on the operator $\tau_p^*\tau_p$ are satisfied (We refer the reader to
Appendix 0
 for the definition of G\^ateaux-differentiability and related notions):
   \item{\bf 1)}
{\it Assumption on the spectral properties of $\tau_{p_0}^*\tau_{p_0}$}
\item{}The operator $e^{-\e \tau_{p_0}^* \tau_{p_0} }$
is trace class for any $\e>0$ and for any vector $X$ at point $p_0$,
there is a neighborhood $\IC_0$ of $p_0$ on the geodesic $p_\kappa=exp_{p_0}
\kappa X$ such that for all $p\in \IC_0$, $e^{-\e \tau_{p_0}^* \tau_{p_0} }$ is
trace class.
 \item{\bf 2)}
 {\it Regularity assumptions }
 \item{}   We shall assume that  the maps
   $p\mapsto \tau_p$ and  $p\mapsto \tau_p^*\tau_p$  are  G\^ateaux
differentiable and that for any $t>0$,
  the
 function $p\mapsto
\hbox{tr} e^{-t\tau_p^*\tau_p}$ is   G\^ateaux
differentiable at point $p_0$.
 \item{}We furhtermore assune that the G\^ateaux-differentials at point $p_0$
 in the direction $X$ of these operators are related as follows:
  $$\delta_X(\hbox{tr} e^{-\e\tau_p^*\tau_p})=
-\e \hbox{tr}
 (\delta_{  X }(\tau_p^*\tau_p) e^{-\e\tau_p^*\tau_p })
  \eqno (2.1)$$
\item{}Moreover, for any vector $X$ at point $p_0$, there  are constants
$C>0$,
$u>0$ and a neighborhood $I_0$ of $p_0$ on the geodesic $p_\kappa
=\hbox{exp}_{p_0} \kappa X$
such that
for any $p\in I_0$:
$$\hbox{tr}e^{-t\tau_p^*\tau_p }\leq C e^{-tu}\eqno (2.2)$$
 and
$$ M_{I_0}(t)\equiv \hbox{sup}_{p\in I_0}\vert\vert\vert \delta_{\bar X(p)}
 (\tau_p^*\tau_p)    e^{-t\tau_p^*
\tau_p}\vert\vert\vert_\infty \eqno (2.3  )$$
  is finite and  a decreasing function  in $t$.
\item{}
Here $\vert\vert\vert \cdot\vert\vert\vert_\infty$ denotes the operator norm
on $\GC$ induced by
$(\cdot, \cdot)_p$,  $\bar X$ is a local vector field  defined in a
neighborhood of $p_0$
by
$\bar X(p_\kappa )= \hbox{exp}_{p_\kappa*}(\kappa X)(X).  $
 \item{}
\item{}The orbit $O_{p_0}$ is called {\it  volume-regularisable  } if
  dim Ker$\tau_p^*\tau_p$ is   constant on some neighborhood of
 $p_0$ on any geodesic contaning $p_0$ and if
   the following assumption is satisfied:
 \item{\bf 3)}
  {\it Assumption on the asymptotic behavior of the heat-kernel traces
 }
\item{}There is  an integer $m>0$ and a family of
 maps
 $p\mapsto  b_j (p),   j\in
\{-J,\cdots ,m-1\} $
 which are
G\^ateaux differentiable in the direction $X$ at point $p_0$  such
 that
$$\hbox{tr} e^{-\e \tau_p^*\tau_p }\simeq_0 \sum_{j=-J}^{m-1} b_j(p)
\e^{j\over m}\eqno (2.4)$$
in a neighborhood $\IC_0$ of $p_0$ on the geodesic $p=exp_{p_0}\kappa
X$,
 and
$$\delta_X \hbox{tr} e^{-\e\tau_p^*\tau_p}\simeq_0
 \sum_{j=-J}^{m-1}\delta_X b_j(p)
  \e^{j\over m}.\eqno (2.5)$$
Furthermore, setting
$$F_p(t)\equiv \hbox{tr} e^{-t\tau_p^*\tau_p}-\sum_{j=-J}^{m-1}
b_j(p)t^{j\over m}$$
for any vector $X$ at point $p_0$, there  is a  constant
$K>0$, and a neighborhood $I_0$ of $p_0$ on the geodesic $\kappa \to
p_\kappa=\hbox{exp}_{p_0} \kappa X$
such that:
$$\hbox{sup}_{p\in I_0 }\Vert \delta_{\bar X(p)}F_p(t)\Vert_\infty \leq
Kt.\eqno (2.5 bis)$$
\item{}A principal bundle as described above with  all its orbits
 volume-preregularisable  (resp. volume- regularisable)
will
be called  preregularisable (resp.  regularisable).      \b
{\bf Remark}: Since   the Riemannian structure on $\PC$ is right invariant and
the one on $\G$ is  $Ad_g$ invariant, the above assumptions do not
 depend
on the point chosen in the orbit for we have
$\tau_{p\cdot g}=R_{g_*}\tau_p {\hbox{Ad}}_g$.\b
 Although most fibre bundles we shall come across are not only preregularisable
but
 also regularisable so that the notion of preregularisabiblity might seem
somewhat
artificial, in   applications (see Appendices A and B), it is often enough
to verify
the conditions required for   preregularisability  in order to prove a
certain
 minimality of the orbits, namely strong minimality, a notion which will be
defined in the following.\s
 Natural examples of regularisable fibre bundles arise in gauge
 field theories (Yang-Mills, string theory). In gauge field theories,
$\PC$ and $\GC$ are modelled on spaces of sections of
vector bundles $\EC$ and $\FC$ based on a compact
 finite dimensional manifold $M$ and the operators
$ \tau_p^* \tau_p$ arise as smooth families of Laplace
 operators on forms. As elliptic operators on a compact
 boundaryless manifold, they have purely discrete spectrum which satisfies
 condition
 1) (see [G] Lemma 1.6.3)
and (2.4) (see [G] Lemma 1.7.4.b)). By classical results
 concerning  one parameter families of heat-kernel operators, they
satisfy (2.1) (see [RS] proposition 6.1) and (2.2) (see
 proof of Theorem 5.1 in [RS]). Since $\delta_XB_p$ is also a partial
differential
operator, by [G] lemma 1.7.7,
$\delta_X \hbox{tr} e^{-\e B_p}$ satisfies (2.5). Assumptions on the
G\^ateaux-differentiability and
assumptions (2.3  ), (2.5 bis) are fulfilled in applications. Indeed,
the parameter $p$ is a geometric object such as a connection, a metric on  $M$
and choosing these objects regular enough (of class $H^k$ for $k$ large enough)
ensures that the maps  $p\mapsto \tau_p$, $p\mapsto \tau_p^*\tau_p$,
 $p\mapsto \hbox{tr} e^{-t\tau_p^*\tau_p}$
etc.. are
 regular enough for they involve these
 geometric quantities and their derivatives, but no derivative of higher
order.\s
{\bf Remark}: In the context of gauge field theories, the underlying
Riemannian structure w.r.to which  the traces (arising in (2.2)-(2.5bis)
are taken are
weak $L^2$  Riemannian structures, the ones that
 also underly the theory of elliptic operators
 on compact manifolds.
In [AP2], we discuss in how far this weak Riemannian structure could be
replaced by
a strong Riemannian structure, in order to set up a link between this geometric
 picture and a stochastic one developped in [AP2].
  \b {\bf Proposition 2.1}:  Let  $O_{p_0}$ be a volume-preregularisable orbit
  such that for any geodesic containing $p_0$, there is a neighborhood of $p_0$
  on this geodesic on which $\tau_p^*\tau_p$ is injective. Then
  \item{1)}
$\hbox{det}_\e(\tau_{p }^*\tau_{p } )$ is well defined for any $\e>0$
and for   $p$ in a neighborhood of $p_0$ on any geodesic of $p_0$.
 \s\m
\item{ 2)} The map
$$p\mapsto \hbox{det}_\e (\tau_p^*\tau_p)$$
 is  G\^ateaux-differentiable at point $p_0$,   the operator
$\int_\e^{+\infty} \delta_X(\tau_p^*\tau_p) e^{-t\tau_p^*\tau_p}dt$
is trace class for any  $p  $ in a neighborhood of $p_0$ on any geodesic of
 $p_0$. For any tangent vector $X$ at point $p_0$, we have:
$$\eqalign{ \delta_X
\hbox{  log det}_\e(\tau_p^*\tau_p)
  &=  \int_\e^\infty  \hbox{tr }   ( \delta_X\tau_p^*\tau_p )
e^{-t\tau_{p_0}^*\tau_{p_0}} dt \cr
&=\hbox{tr} \int_\e^{+\infty} (\delta_X \tau_p^*\tau_p)
 e^{-t\tau_p^*\tau_p}dt\cr}\eqno (2.6 a)  $$
  \item{ 3)}If the orbit $O_{p_0}$ is  moreover volume-regularisable,
the map $p\mapsto \hbox{det}_{reg} (\tau_p^*\tau_p)$ is  G\^ateaux
differentiable in all directions at point $p_0$, and with the notations of
(2.4)
   $ \delta_X[\hbox{    log det  }_\e(\tau_p^*\tau_p)-
\sum_{j=-J }^{-1} {m\over j}
 b_j\e^{j\over m}
-  b_0 \hbox{log} \e ] $
 converges when $\e \to 0$  and we have
 $$\eqalign{&\lim_{\e \to 0}
\delta_X [\hbox{  log det  }_\e(\tau_p^*\tau_p) -
\sum_{j=-J }^{-1} {m \over j} b_j \e^{j\over m}-
b_0 \hbox{log} \e ]
 = \delta_X\hbox{  log  det}_{reg}\tau_p^*\tau_p\cr
&= -\sum_{j=-J,j\neq 0}^{m-1} {m\over j}\delta_Xb_j-\int_1^\infty
\delta_X\left(\hbox{tr}
{e^{-t\tau_p^*\tau_p}\over t}\right) dt-\int_0^1
{\delta_X F_p(t) \over t} dt \cr} \eqno (2.6 b)
 $$
 \b
 \m {\bf Proof  }:We set $B_p=\tau_p^*\tau_p$.
\item{1) }By the first assumption for volume-preregularisable orbits,
we know that
$e^{-\e B_p}$ is trace class so that by lemma 1.1 so is
  $A_\e^p\equiv \hbox{log} h_\e(B_p)$. Hence
$\hbox{det}_\e(B_p)=e^{\hbox{tr} A_\e^p}$ is well defined.
 \item{2)}  Let us show the first equality in (2.6 a). The orbit $O_{p_0}$
being preregularisable,
  assumption  2) for volume-preregularisability  yields that
for any $p\in I_0$ and any
 $t>\e>0$
  $$  \vert \hbox{tr}(\delta_{\bar X(p )}B_p
e^{-tB_{p }})   \vert \leq C M_{ I_0}({t\over 2}) e^{-{t\over 2} u}. $$
Here, we have used the fact that
$\vert \hbox{tr}(UV)\vert\leq \vert\vert \vert U\vert\vert \vert \hbox{tr} V
 \vert$ for any bounded operator $U$ and any
trace class operator
$V$ applied to
 $U= \delta_{\bar X(p)} B_p e^{-{t\over 2}B_{p_0}}$ and $V=e^{-{t\over 2}
B_p}$.
Hence, by Lebesgue dominated convergence theorem,
the map $p\mapsto \int_\e^\infty t^{-1}\hbox{tr} e^{-tB_p}dt$
is G\^ateaux-differentiable in the direction $X$ at point  $p_0$ and
$$\eqalign{\delta_X \int_\e^\infty t^{-1}\hbox{tr} e^{-tB_p}dt&=
 \int_\e^\infty t^{-1}\delta_X\hbox{tr} e^{-tB_p}dt\cr
&= -\int_\e^\infty  \hbox{tr}((\delta_XB_p) e^{-tB_{p_0}})dt\cr}$$
using (2.2).
 Using the fact that
 $\hbox{log det}_\e(B_p)=
-\int_\e^{+\infty}t^{-1} \hbox{tr} e^{-tB_p}dt$
then yields the first equality in (2.6 a).
\item{}The second  equality in (2.6 a) and the fact that   we can swap the
trace  and the integral
follow  from the estimate:
$$ \vert\vert \vert  \delta_X B_p e^{-tB_{p_0}}  \vert\vert \vert_1
 \leq \vert\vert \vert \delta_XB_p e^{{-\e \over 2}B_{p_0}}\vert\vert
\vert_\infty \Vert \vert
  e^{-\half tB_{p_0}}  \Vert \vert_1 \leq C\vert\vert \vert \delta_XB_p
e^{{-\e \over 2}B_{p_0}}\vert\vert \vert_\infty \Vert \vert_\infty
  e^{-tu} \eqno(*)$$
valid for $t\geq \e$, using assumption (2.2). We finally  obtain by
dominated convergence:
$$ \hbox{tr} \int_\e^{+\infty} \delta_X B_p e^{-tB_{p_0}} dt
  =\int_\e^{+\infty} \hbox{tr} \delta_X B_p e^{-tB_{p_0}} dt. $$
 \item{3)}
  Let us first check that the map $p\mapsto \hbox{det}_{reg} B_p$ is
G\^ateaux differentiable at point $p_0$ in the direction $X$.
By (1.6), we have
$$ \hbox{  log  det }_{reg} B_p=-\sum_{j=-J, j\neq 0 }^{m-1} {b_j(p)\over j}
-\int_1^\infty \hbox{tr} {e^{-tB_p}\over t} dt- \int_0^1 {F_p(t)\over t} dt$$
The first    term  on the r.h.s. is  G\^ateaux differentiable in the
direction $X$ by the
assumption on the maps $p\mapsto b_j(p)$.
The second term  on the r.h.s. is
 G\^ateaux differentiable by
the result (applied to $\e=1$  of part 2) of this lemma which
 tells us that $p\mapsto \hbox{det}_\e(B_p)$
 is  G\^ateaux differentiable.
The G\^ateaux differentiability of the   last
term follows from the local uniform upper bound (2.5 bis).\m  \item{}We now
check (2.6 b). The map $p\mapsto \hbox{log det}_\e(B_p)-
\sum_{j=-J , j\neq 0}^{ m-1} {mb_j\over j} \e^{j \over m}-b_0\hbox{log} \e$
is  G\^ateaux differentiable  in the direction $X$ and  we can write
  $$\eqalign{ & \delta_X\hbox{   ( log det}_\e(B_p)-\sum_{j=-J,j\neq 0}^{m-1}
{mb_j \e^{j\over m}\over j}-b_0 \hbox{log} \e) \cr
&=\delta_X (-\int_\e^{\infty} \hbox{tr}
{e^{-tB_p}\over t}dt -\sum_{j=-J,j\neq 0}^{m-1} {mb_j \e^{j\over m} \over j}
-b_0 \hbox{log} \e) \cr
 &=\delta_X  \left( -\sum_{j=-J,j\neq 0  }^{ m-1}
 m{b_j\over j}- \int_1^\infty \hbox{tr} {e^{-tB_p}\over t} dt-
\int_\e^1 {F_p(t)\over t} dt \right) \quad \hbox{as in (1.9)}\cr
&=    -\sum_{j=-J,j\neq 0 }^{m -1}
 \delta_Xb_j{m\over j} -
 \int_1^\infty \delta_X\hbox{  tr} {e^{-tB_p}\over t}  dt-
 \int_\e^1 \delta_X{F_p(t)\over t} dt
\cr}$$
which tends to
$\delta_X\hbox{  log det}_{reg} B_p $ by (1.6)  and dominated convergence.
Here we have used  the results of   point  2) of the proposition applied to
$\e=1$  to write
$$\delta_X\int_1^\infty \hbox{tr}{ e^{-tB_p}\over t}dt =
 \int_1^\infty  \delta_X\hbox{  tr} e^{-tB_p}  dt $$
and (2.5 bis) to write     $\delta_X\int_\e^1 {F_p(t) \over t}dt=
\int_\e^1 {\delta_X F_p(t) \over t} dt$.
 \kasten \b
{\bf Remark} These results extend to the case when instead of assuming that
$\tau_p^*\tau_p$ is injective locally around $p_0$,
  one considers orbits of an action at points $p_0$ for which the
dimension
of the kernel
of $\tau_p$ is constant on some neighborhood of $p_0$ on each geodesic
starting at point $p_0$. For this, one should replace
$\hbox{det}_\e\tau_p^*\tau_p$ and $\hbox{det}_{reg}\tau_p^*\tau_p$ by
$\hbox{det}_\e^\prime\tau_p^*\tau_p$ and
$\hbox{det}_{reg}^\prime\tau_p^*\tau_p$.
This extension is useful for applications (see Appendix A.).
\b
A direct generalisation of the notion of volume for volume-preregularisable
or regularisable orbits would give infinite quantities. But for
volume-preregularisable
or regularisable orbits, one can define a  notion of
  preregularised
or regularised volume, which justifies a posteriori the term
"volume-preregularisable
or volume-regularisable orbits" for these orbits.
 Since
$\tau_{p\cdot g}=R_{g*} \tau_p \hbox{Ad}_g$ and since  the metric on $\G$
is $\hbox{Ad}_g$
and that on $\PC$ right  invariant,  for any $\e>0$, we have
$\hbox{det}_\e(\tau_{p\cdot g}^* \tau_{p\cdot g})=
\hbox{det}_\e(\tau_p^*\tau_p)$
so that it makes sense to set the following definitions:
\m
 {\bf Definition:}
\item{1)} Let $O_p$ be   a  volume-preregularisable orbit, then
$$\hbox{vol}_\e(O_p)\equiv\sqrt {\hbox{det}_\e^\prime(\tau_p^*\tau_p)}$$
defines a one parameter family of (heat-kernel) preregularised
volumes of $O_p$.
\item{2)} Let $O_p$ be   a   volume-regularisable orbit, then
$$\hbox{vol}_{reg}(O_p)=
\sqrt {\hbox{det}_{reg}^\prime(\tau_p^*\tau_p)}$$
defines the heat-kernel regularised volume of $O_p$.
\item{3)} Let $O_p$ be   a   volume-regularisable orbit, then
$$\hbox{Vol}_{reg}(O_p)=
\sqrt {\hbox{Det}_{reg}^\prime(\tau_p^*\tau_p)}$$
defines the zeta function volume-regularised volume of $O_p$.
\m From lemma 1.4 follows that
$$\hbox{Vol}_{reg}(O_p)=e^{ -\half\gamma b_0^\prime(p)} \hbox{vol}_{reg}
(O_p)$$
where $\gamma$ is the Euler constant and $b_0^\prime(p)=b_0(p)-\hbox{dim Ker}
(\tau_p^*\tau_p)$ is the coefficient
arising from the heat-kernel asymptotic expansion of $\tau_p^*\tau_p$
given by (2.4).\s
{\bf Remarks}: \item{1)}In finite dimension, dim$H=d$ and $\tau_p$ is
injective, we have by (1.13):
$$\eqalign{\lim_{\e \to 0}( \e^{-  d\over 2}\hbox{vol}_\e(O_p))&=
\hbox{vol}_{reg}(O_p)=  e^{d\gamma \over 2} \hbox{Vol}_{reg} (O_p)\cr
&=
 e^{d\gamma\over 2} \vert\hbox{det}\tau_p\vert\int_G \vert \hbox{det}
\hbox{Ad}_g d\mu(g)\vert = e^{d\gamma\over 2} {\hbox{vol}(O_p)\over
\hbox{vol}G}
 \cr}
\eqno (2.7)$$
where $\mu$ is the volume measure and  vol$(O_p)$
  is  the ordinary volume of
 the fibre
 $O_p$.
\item{2)}
If the coefficients $b_j(p)$ arising in the heat-kernel expansion of
 $\tau_p^*\tau_p$ are independent of $p$, we have for two regularisable
orbits $o_{p_0}$ and $O_{p_1}$:
$$ \lim_{\e \to 0}{\hbox{vol}_\e(O_{p_0}) \over \hbox{vol}_\e(O_{p_1})} =
 {\hbox{vol}_{reg} (O_{p_0}) \over \hbox{vol}_{reg} (O_{p_1})}$$
 \m
{\bf Proposition 2.2}: The heat-kernel (pre)-regularised and zeta function
regularised
volume
  of a volume-(pre)regularisable orbits $O_p$ is  G\^ateaux-differentiable
at the point $p$.\m{\bf Proof}: It follows from Proposition 2.1.\kasten\b
Let us now introduce a notion of extremality of  orbits  which
 generalises the corresponding finite dimensional notion [H].\m
{\bf Definition}: A    {\it   strongly extremal orbit} is a
volume-preregularisable
   orbit, the heat-kernel preregularised volume of which is
extremal, i.e $O_p$ is strongly extremal if  $\delta_X \hbox{vol}_\e(O_p)=0$
for
  any horizontal vector $X$ at point $p  $ and any $\e>0$.
\s A  heat-kernel (resp. zeta function)   {\it extremal orbit} of a
volume-preregularisable bundle
is an orbit, the heat-kernel (respectively zeta function) regularised
volume of which is
extremal, i.e $\delta_X \hbox{vol}_{reg}(O_p)=0$ for any horizontal vector
$X$ at point
$p $.
\s
Notice that whenever $b_0 $ does not depend on $p$,
the zeta-function regularised volume of an extremal orbit is also heat-kernel
extremal. From  (2.7) also follows that this notion generalises the finite
dimensional notion of extremality of the volume of the fibre.
\m
 \b  \b {\bf  III. Minimal regularizable orbits as orbits with extremal
regularized volume }
\m
We shall consider a  preregularisable principal fibre bundle $\PC\to \PC/\G$.
  By assumption, the bundle  is equipped with a Riemannian
connection given by a family of horizontal spaces $H_p, p\in \PC$ such that
 $$T_p\PC=H_p\oplus V_p$$
where $V_p $
is the tangent space to the orbit at point $p$ and the sum is an orthogonal
one.
  \m
For a horizontal vector   $X$ at point $p$, we define the shape operator
$$\eqalign{ \H_{ X }:V_p&\to V_p\cr
Y&\mapsto  -(\nabla_Y\bar X)^v(p)\cr}$$
where the subscript $v$ denotes   the orthogonal projection onto $V_p$
and
$\bar X$ is a horizontal field with value $X$ at $p$. Similarly, we define
the second fundamental form:
$$\eqalign{ S^p:V_p\times V_p&\to H_p\cr
(Y, Y^\prime) &\mapsto (\nabla_{\bar Y}\bar Y^\prime)^h(p)\cr}$$
where $\bar Y$, $\bar Y^\prime$ are vertical vector  fields such that
$\bar Y(p)=Y$, $\bar Y^\prime(p)=Y^\prime$. These definitions are
 independent of the choice of the extensions of $X$,$Y$ and $Y^\prime$.
\s An easy computation  shows that the shape operator and the second
fundamental form are related as follows:
$$<\H_{ X }(Y), Y^\prime>_p=<S^p(Y, Y^\prime),   X >_p \eqno (3.1)$$
 Note that this  explicitely shows that $\HC_{  X }$ only depends on $X$ and
not on the
  extension
 $  \bar X$ of $X$.  Since $S^p$ is symmetric, so is $\HC_X$.\s
As in the finite dimensional case, one can define the notion of
totally geodesic orbit, an orbit $O_p$ being totally geodesic whenever
the second fundamental form $S^p$ vanishes.
 \b {\bf Definition}:  The orbit $O_p$  of a point $p\in \PC$ will be called
{\it heat-kernel preregularisable } or for short {\it preregularisable} if
for   any horizontal  vector   $ X$ at $p$, $\f \e>0$,
$$\H_{ X }^\e\equiv e^{-\half\e \tau_{p}\tau_{p}^*}\H_{  X }
e^{-\half \e \tau_{p}
\tau_{p}^*}\eqno (3.2 )$$
 is trace class. A preregularisable orbit $O_p$ will
be called {\it   strongly
minimal } if moreover for any  $q\in O_p$ and  $X$ a horizontal vector at
point $q$,
 tr$\HC_{  X }^\e=0$ $\f \e>0$. \m
{\bf Remarks}:     \item{1)} The preregularisability of the orbits ( namely
$\HC_X^\e$ trace class)
 is automatically satisfied if the manifold $\PC$ is equipped with a strong
smooth Riemannian structure, since in that case the second fundamental form
is a bounded bilinear form
and its weighted trace is well defined
(see also [AP2] where this is discussed in further details).
\item{2)} Since on a preregularisable bundle, the Riemannian structure on
$\PC$ is right invariant and
the one on $\G$ is  $Ad_g$ invariant, the notion of
 (pre) regularizability and (strong) minimality   of the orbit does not
depend on  the
 orbit chosen on the orbit. Indeed, let $p$ be a point, $g\in \G$ and $X$ a
horizontal vector at point $p$.
Let $\bar X$ be a right invariant horizontal vector field coinciding with
$X$ at $p$.
Since
$\tau_{p\cdot g}=R_{g_*}\tau_p {\hbox{Ad}}_g$ and $
\HC_{\bar X(p\cdot g)}=R_{g*}\HC_{\bar X(p)} R_{g*}^{-1}$, the vector field
$\bar X$ being right invariant,
  we have
$\tau_{p\cdot g}\tau_{p\cdot g}^*=R_{g*} \tau_p
\tau_p^* R_{g*}^{-1}$.
 Hence  $\H_{\bar X(p\cdot g)}^\e=R_{g*} \H_{\bar X(p)}^\e R_{g*}^{-1}$ is
trace class w.r. to
$<\cdot, \cdot >_{p\cdot g}$ whenever $\H_{\bar X(p)}^\e$ is trace class
 w. r. to
$<\cdot ,\cdot >_p$ and
$\hbox{tr} \HC^\e_{\bar X(p\cdot g)}=\hbox{tr} \HC_{\bar X(p)}^\e$.
 \item  {3)} Notice that if $\HC_X$ is trace class, as in the  finite
dimensional case,
strong minimality implies that $\hbox{tr} \HC_X=0$ and
hence ordinary minimality. The fact that strong minimality implies
minimality in the finite dimensional case motivates the choice of the
adjective "strong".
\item{4)} This preregularised shape operator $\HC^\e_X$ and the second
fundamental
form are related as follows:
$$<\HC_X^\e(Y),Y^\prime>_p= <S^p(e^{-\half\e \tau_p\tau_p^*}Y,
e^{-\half\e\tau_p \tau_p^*}Y^\prime), X>_p$$
  Since $\tau_p \tau_p^*$ is an isomorphism of the tangent space to the
 fibre
$T_pO_p$, $\HC_X^\e$ vanishes  whenever the second fundamental form
vanishes and an
 orbit is  totally geodesic whenever this regularised shape operator
 vanishes on the orbit for some  $\e>0$.
 \b {\bf Definition}: A preregularizable orbit $O_p$ will be said to be
{\it regularisable}
 if furthermore,
 the one parameter family
 $\HC_{ X }^\e, \e\in ]0,1]$ admits a regularized limit-trace
$\hbox{tr}_{reg} \HC_{  X }$.
  \s  For a preregularizable orbit $O_p$ such that for any $\e>0$,
$X\to \hbox{tr} \HC_X^\e$ is a bounded linear form on  $T_p\PC$ for
 the norm induced by $<\cdot, \cdot>_p$,
  by Riesz theorem we can define the
 {\it preregularised mean    curvature
vector}
   $  \SC_\e$ in the closure $H_p$ of $T_p\PC$ for this norm   by
 the relation
 $$ <\SC_\e(p),   X  >_p= \hbox{tr}  \H_{  X }^\e\eqno (3.3)$$
In the same way, for a regularisable
 orbit $O_p$ such that
$X\to \hbox{tr}_{reg} \HC_X$ is a bounded linear form on $T_p\PC$,
by Riesz theorem, we can
define
 the {\it regularized mean   curvature
 vector} $\SC_{reg}(p)$ in $H_p$
   by
 the relation
 $$ <\SC_{reg}(p),   X  >_p= \hbox{tr}_{reg} \H_{  X }\eqno (3.4)$$
for any horizontal vector $X$ at point $p$.
Of course, if the Riemannian structure is strong, both $\SC_\e(p)$ and
$\SC_{reg}(p)$ lie in $T_p\PC$.
\s{\bf Remark}: In the finite dimensional case,
we have
$\SC_0(p)=\SC_{reg}(p)=\hbox{tr}S^p$ where $S^p$ is the second fundamental
form.
In the infinite dimensional case, the family of preregularised mean
 principal curvature vectors
 $\SC_\e(p)$ coincides with a family of preregularised traces and the
regularised mean principal curvature vector
$\SC_{reg}(p)$ with the
 regularised
trace of the second fundamental form $S^p$.
  \s {\bf Definition}: A regularisable orbit $O_p$  will be called {\it
heat-kernel
  minimal}
  if $\hbox{tr}_{reg}\H_{  X }  =0 $ for any horizontal vecotr at point $X$.
\s {\bf Remarks}:
\item{1)}Since on a preregularisable bundle, the Riemannian structure on
$\PC$ is right invariant and
the one on $\G$ is  $Ad_g$ invariant,
  regularizability and  minimality   of the orbit $O_p$ does not depend on the
point $p$ chosen on the orbit.
As before, we have $\hbox{tr}_{reg} \HC_{\bar X(p\cdot g)}=\hbox{tr}_{reg}
 \HC_{\bar X(p)} $.
\item{2)} Here again, in the finite dimensional case, the one parameter
 family
$\HC_X^\e$ admits a regularised limit trace given by the ordinary trace
$\hbox{tr}_{reg}\HC_X=
\hbox{tr} \HC_X$ and heat-kernel minimality is equivalent to
the finite dimensional notion of minimality.
 \m Note that a strongly minimal preregularisable orbit $O_p$ is
regularisable and
minimal since setting
$A_\e\equiv \HC_{ X }^\e$ in (1.1), we have $a_j=0, \f j$ and hence
$\hbox{tr}_{reg} \HC_{ X }=0$.
\s
 \b\b The regularisation of
 the mean principal curvature vector for
orbits of  group actions in the infinite dimensional case has been
 discussed in the literature before.
 King and Terng in [KT] introduced a notion of regularisability and
minimality for submanifolds  of path spaces    using zeta-function
regularisation methods. They in particular show zeta function
 regularisability and minimality for
 the orbits of the coadjoint action of a  (based) loop group on a space of
loops in
 the corresponding Lie algebra.
We  shall show later on that these orbits within this framework
 are   regularisable and  strongly minimal (hence   minimal).
\s We now introduce a notion of  zeta function
regularisability  which is a slight variation of the one introduced by
 Maeda, Rosenberg and Tondeur in [MRT1] (see also [MRT2]) in the case of
   orbits of  the gauge action in Yang-Mills theory.
This modification is natural in our context as we shall see later on.
 \s {\bf Definition}: The orbit $O_p$ of a point $p$ is {\it
zeta function regularisable } whenever
 $$\lim_{s\to 1}
-\half\left[ \Gamma(s)^{-1} \int_0^\infty t^{s-1} \sum_{\l_n \neq
0}e^{-t\l_n^p }
\delta_X \l_n^p dt + (s-1)^{-1} \delta_Xb_0^\prime (p)dt\right]$$
exists for any horizontal field $X$ at point $p$
and the limit shall be denoted by $Tr_{reg} \HC_X$ so that
$$ Tr_{reg} \HC_X=-\half \lim_{s\to 1}
\left[ \Gamma(s)^{-1} \int_0^\infty t^{s-1} \sum_{\l_n\neq 0} e^{-t\l_n^p}
\delta_X \l_n^pdt + (s-1)^{-1}\delta_Xb_0^\prime (p) \right]\eqno (3.5)$$
 Since on a regularisable principal fibre bundle, the Riemannian structure
on $\PC$ is $\G$ invariant and that on $\G$ is Ad$\G$ invariant
  the orbit  $O_p$ is zeta function regularisable whenever the above
 holds for any  $q\in O_p$. \m
 \s {\bf Definition}: A zeta function regularisable orbit $O_p$  will be
called {\it zeta
 function
  minimal}
  if $\hbox{Tr}_{reg}\H_{  X }  =0 $ for any horizontal vecotr at point $X$.

  \m
{\bf Remark}: This definition coincides with  that of [MRT
  Proposition 5.9] whenever $\delta_X b_0(p)$ is zero. This notion of
regularisability is of course less restrictive than
that of [MRT] and we shall see that on a regularisable fibre bundle, there
is no obstruction to zeta function regularisability of the orbits.\s
 \m
\b Let us introduce some notations. Let $\PC\to \PC/\G $ be a
preregularisable principal
fibre bundle and
let $(T_n^p)_{n\in \N}$  be  a   set
 of eigenvectors of $\tau_p^*\tau_p$ in $\GC$ corresponding to the
 eigenvalues
$(\l_n^p)_{ n\in \N}$  counted with multiplicity and in
 increasing order.
Let $p_0$ be a fixed point in $\PC$ and let  $\IC^p_{p_0}$ be  the
isometry from
$(\GC, (\cdot, \cdot)_{p_0})$ into $(\GC, (\cdot, \cdot)_p)$
 which takes the orthonormal set  $(T_n^{p_0})_n$  of eigenvectors of
$\tau_{p_0}^*\tau_{p_0}$ to the orthonormal set of eigenvectors
$(T_n^p)_n$ of $\tau_p^*\tau_p$.
Notice that $\IC_{p_0}^{p_0}=I$.
\b {\bf Lemma 3.1}: Let $\PC \to \PC/\GC$ be a
preregularisable principal fibre bundle.
Let $p_0\in \PC$ be a point at which the map $p\mapsto \IC_{p_0}^pu$
is G\^ateaux-differentiable for any $u\in \GC$. Let $X $ be a horizontal
vector
 at   $p_0 $.  We shall consider eigenvalues $\l_n^p$ that correspond to
eigenvectors that do not belong to
 $\IC_{p_0}^p \hbox{Ker} \tau_{p_0}^*\tau_{p_o}$.\s
  \item{1)}  The maps $p\to \l_n^p$ are
G\^ateaux-differentiable in the direction $X$ at point $p_0$,
$$\delta_X \l_n^p=
 (\delta_X (\tau_p^*\tau_p) T_n^{p_0}, T_n^{p_0})_{p_0}$$
and
$$\delta_X \hbox{log}h_\e(\l_n^p)=
  \int_\e^{+\infty}(\delta_X (\tau_p^*\tau_p) e^{-t \tau_{p_0}^*
\tau_{p_0}}T_n^{p_0}, T_n^{p_0})_{p_0}dt.$$
    \item{2)}Furthermore,we have
 $$
 -  <\H_{ X }^\e\tilde U_n^{p }, \tilde U_n^{p }>_{p_0}+
  e^{-\e \l_n^{p_0}}(\delta_X\IC^p_{p_0}  T_n^{p_0},
 T_n^{p_0})_{p_0 }
  =\half \delta_X\hbox{  log } h_\e (\l_n^p)  \eqno(3.6a)$$
 where we have set
$\tilde U_n^p=\Vert  \tau_p T_n^p\Vert^{-1}\tau_p T_n^p$.
\item{3)} If the Riemannian structure on $\GC$ is fixed
 (independent of $p$), then $\delta_X \IC_{p_0}^p $ is antisymmetric
and
$$\eqalign{&\half   \int_\e^{+\infty} ( \delta_X(\tau_p^*\tau_p)
e^{-t\tau_p^*\tau_p
 } T_n^p,
T_n^p)_{p_0} dt=-  <\H_{ X }^\e\tilde U_n^{p }, \tilde U_n^{p }>_{p_0}
  \cr
 &=\half \delta_X\hbox{  log } h_\e (\l_n^p)=
\half \l_n^{{p_0}^{-1}} \delta_X \l_n^p e^{-\e \l_n^{p_0}} \cr}\eqno (3.6  b)$$
     \m {\bf Proof}:
As before, we shall set $B_p=\tau_p^*\tau_p$.
 Since $p_0$ is fixed, we drop the index $p_0$ in $\IC_{p_0}^p$ and
denote this isometry by $\IC^p$. Notice that $\IC^{p_0}=I$.
As before, we denote by
$(T_n^p)_{n\in \N}$ the orthonormal set of eigenvectors of
 $\tau_p^*\tau_p$ which correspond  to the eigenvalues $(\l_n^p)_{n\in \N}$
in increasing order and counted with multiplicity.
We shall set $\tilde T_n^p= \tau_p T_n^p$, $\bar T_n^p= \tau_p T_n^{p_0}$.
  \item{1)} Using the relations
 $(\IC^p\cdot, \IC^p\cdot)_p=(\cdot, \cdot)_{p_0}$,
  $\IC^p(T_n^{p_0})=T_n^{p }$,
${\IC^p}^*\IC^p=I$, we can write
$\l_n^p=(B_pT_n^p, T_n^p)_p=(  B_p \IC^p T_n^{p_0},
 \IC^pT_n^{p_0})_{p_0}$
 and the map $p\mapsto \l_n^p$ is   G\^ateaux
 differentiable in all directions at point $p_0$
 since
$p\mapsto B_p$, $p\mapsto \IC^p$ are   G\^ateaux-differentiable by
assumption on the bundle. Furthermore
$$\eqalign{ \delta_X(B_pT_n^p, T_n^p)_p&=
\delta_X( {\IC^p}^* B_p {\IC^p}  T_n^{p_0}, T_n^{p_0})_{p_0}\cr
&=((\delta_X B_p) T_n^{p_0}, T_n^{p_0})_{p_0} +(\delta_X
 ({\IC^p}^* )
B_{p_0}
T_n^{p_0 }, T_n^{p_0})_{p_0} +\cr
&+ ( {\IC^{p_0}}^* B_{p_0} (\delta_X{\IC^p} ) T_n^{p_0}, T_n^{p_0})_{p_0}
\cr
 &=((\delta_X B_p) T_n^{p_0}, T_n^{p_0})_p+
\l_n^{p_0}([{\IC^{p_0}}^*\delta_X ({\IC^p}  )  +
  (\delta_X{\IC^p}^*)  {\IC^{p_0} } ] T_n^{p_0 } , T_n^{p_0 })_{p_0 }\cr}$$
Since
 ${\IC^p}^*\IC^p=I$, we have
 $ \delta_X{\IC^p}^* \IC^{p_0} +{\IC^{p_0}}^*\delta_X {\IC^{p }}  =0$ so that
finally $\l_n^p$ is  G\^ateaux-differentiable and
 $\delta_X\l_n^p=((\delta_X B_p) T_n^{p_0}, T_n^{p_0})_{p_0}$.
\item{}
Using the local uniform estimate (2.3  ), and with the same notations,
 we have for $t>\e$:
 $$ \Vert(\delta_{\bar X(p)}
(B_p)e^{-tB_{p_0}}T_n^{p_0}, T_n^{p_0})_{p_0}\Vert\leq
M_{I_0}(\half t)e^{-\half t\l_n^{p_0}}$$
 so that
the map
$p\mapsto \hbox{log}h_\e(\l_n^p)$ is G\^ateaux-differentiable at point
$p_0$ in the direction $X$
and
$$\eqalign{\delta_X \hbox{log} h_\e(\l_n^p)&=
 -\delta_X \int_\e^\infty t^{-1}(e^{-tB_p} T_n^p, T_n^p)dt\cr
&= (\int_\e^{+\infty} \delta_X(B_p) e^{-t
B_{p_0}}T_n^{p_0}, T_n^{p_0})_{p_0} dt\cr}$$
\item{2)}   \s By definition of $h_\e$ we have:
$$\eqalign{ \delta_X\hbox{ log } h_\e (\l_n^p) &=
(\hbox{log}h_\e)^\prime(\l_n^p)\delta_X \l_n^p  \cr
&= (\l_n^{p_0})^{-1} e^{-\e \l_n^{p_0}} \delta_X \l_n^p \cr}$$
But, with the notations of Appendix 0:
$$\eqalign{ \delta_X \l_n^p &=\delta_X
<\tilde T_n^p, \tilde T_n^p>_p
=2 <\delta_X(\tau_p \IC^p)T_n^{p_0} , \bar T_n^{p_0} >_{p_0}  \cr
 &= 2< \delta_X   \bar T_n^{p }  , \bar T_n^{p_0}>_{p_0}
+2< \tau_p\delta_X \IC^p T_n^{p_0}, \bar T_n^{p_0}>_{p_0}\cr
&=-2 <\nabla_{ \bar T_n^{p_0}} \bar X, \bar  T_n^{p_0} >_{p_0} +2
 < \tau_p\delta_X \IC^p T_n^{p_0},  \bar  T_n^{p_0}>_{p_0}\cr
&=-2 <\nabla_{ {\tilde  T_n}^{ p_0}} \bar X, {\tilde   T_n}^{p_0} >_{p_0}
 +2 < \tau_p\delta_X \IC^p T_n^{p_0}, \bar T_n^{p_0}>_{p_0}\cr
 &= -2\l_n^{p_0}<\HC_X {\tilde  U_n}^{p_0}  ,{\tilde  U_n}^{p_0}>_{p_0}+
2\l_n^{p_0}(  \delta_X \IC^p T_n^{p_0},   T_n^{p_0})_{p_0}  \cr} $$
 where for the third equality,we have used the fact that,
 $\bar X$ being right invariant, we have $[\bar T_n^p, \bar X]=0$.
\s Hence we have:
$$\eqalign{\delta_X\hbox{  log} h_\e (\l_n^p) &=
- 2e^{-\e \l_n^{p_0}}
<\H_{\bar X(p_0)}\bar U_n^{p_0}, \bar U_n^{p_0}>_{p_0}+\cr+
&2
  e^{-\e \l_n^{p_0}}( \delta_X \IC^p
 T_n^{p_0},   T_n^{p_0})_{p_0}  \cr}$$
which yields    2).
\item{3)} On one hand, since the scalar product on the Lie algebra is fixed,
we have $\delta_X\IC^{p*}\subset(\delta_X \IC^p)^*$ (see Appendix 0).
On the other hand,  since $\IC^{p*}\IC^p=I$, we have
$-\delta_X \IC^p\subset \delta_X\IC^{p*} $ (see Appendix 0) so that the
second term in
 the l.h.s of (3.6a) vanishes.
\kasten\b
 \b
{\bf Definition}:
We shall call an orbit $O_{p_0}$ of a preregularised bundle an
 {\it orbit of   type $(\TC)$} whenever the following conditions are satisfied:
\item{\bf 1)}The map $p\mapsto \IC_{p_0}^p$ is G\^ateaux-differentiable
 at  point $p_0 $.
\item{\bf 2)} The operator
$\delta_X \IC_{p_0}^p e^{-\e \tau_{p_0}^*\tau_{p_0}}$
is trace class for any $p_0 \in \PC$ and $\e>0$.
\item{\bf 3)} For any $p\in \PC$,
 $\hbox{tr}\IC_{p_0}^p e^{-\e \tau_{p_0}^*\tau_{p_0}}$ is
G\^ateaux-differentiable at point $p_0\in \PC$
and
$$\delta_X \hbox{tr}
(\IC_{p_0}^p e^{-\e \tau_{p_0}^*\tau_{p_0}})=\hbox{tr}
(\delta_X \IC_{p_0}^p e^{-\e \tau_{p_0}^*\tau_{p_0}}) $$
 \b Whenever the Riemannian structure on $\GC$ is independent of $p$,
any orbit satisfying condition 1) is of type $(\TC)$, for in that case
 the traces involved in 2) and 3)  vanish, $\delta_X \IC_{p_0}^p$
 being an antisymetric operator.
\m Let us interpret the trace tr$(\delta_X\IC_{p_0}^p
 e^{-\e \tau_{p_0}^*\tau_{p_0}})$
as a variation of a relative volume.
Since  the map
 $\kappa\to (\IC^{p_\kappa}_{p_0}\cdot, \IC^{p_\kappa}_{p_0}\cdot)_{p_0}$ is
  differentiable at $\kappa=0$, it is continuous at this point. Hence
for a family  of points $p_\kappa=
\hbox{exp}_{p_0}  (\kappa X)$
 on the geodesic at point $p_0$ generated by $X$, there is a constant
 $\eta >0$ such that for $\alpha $ small enough,
 $\hbox{sup}_{\kappa\in [0, \alpha]} \Vert \IC_{p_0}^{p_0}-
\IC^{p_\kappa}_{p_0}\Vert
 \leq \eta$ (where the operator norm is taken w.r.to $\Vert\cdot \Vert_{p_0}$)
  so that for $\kappa $ small enough
$  \vert ( T_n^{p_\kappa},
T_n^{p_\kappa})_{p_0}\l_n^{p_0} \vert\geq (1-2\eta) \l_n^{p_0} $. Thus
  $\sum_n \hbox{log} h_\e  [( T_n^{p_\kappa},
T_n^{p_\kappa})_{p_0}\l_n^{p_0}]$ converges, since $h_\e$ is non decreasing.
We can therefore  define   a notion of { \it  preregularised
relative volume } of the orbit
 $O_{p_\kappa}$ with respect to $O_{p_0}$:
$$\hbox{vol}_\e^{p_0} (O_{p_\kappa})=\prod_{\l_n\neq 0}
\sqrt { h_\e(  \l_n^{p_0} ( T_n^{p_\kappa}, T_n^{p_\kappa})_{p_0} )}$$
Notice that it coincides with $\hbox{Vol}_\e(O_{p_0})$ for $\kappa=0$.
The fact that $\vert\vert \vert \delta_X \IC_{p_0}^p\vert \vert \vert$ is
locally
 bounded on a geodesic starting at point $p_0$
implies that the
 map $\kappa \to \hbox{log vol}_\e^{p_0} (O_{p_\kappa})$ is
 differentiable at point $\kappa=0$
and an easy
computation yields:
$$\delta_X \hbox{log vol}_\e^{p_0} (O_p)=
\hbox{tr}^\prime[ \delta_X\IC_{p_0}^p e^{-\e \tau_{p_0}^*\tau_{p_0}} ]
$$
 where $tr^\prime$ means that we have restricted to the orthogonal of
Ker$\tau_p^*\tau_p$.\s
 In the next proposition, we investigate the relation between the shape
operator and
 the variation of the volume of the orbit.
\m
{\bf Proposition 3.2 }: Let $\PC\to \PC/\G$ be a
 preregularisable principal
 fibre bundle. Then\m
\item{1)} any orbit   of type $(\TC)$ is preregularisable.
\item{} More precisely, if $O_{p_0}$ is  an
orbit   of type $(\TC)$, for any horizontal vector $X$ at point $p_0$,
the operator $\HC_X^\e$ is trace class,
   the maps   $p\mapsto \hbox{vol}_\e^{p_0}(O_p)$ and
  $p\mapsto \hbox{vol}_\e(O_p)$ are G\^ateaux
 differentiable in the direction $X$ at point $p_0$ and
$$\hbox{tr} \H_{  X }^\e-\delta_X \hbox{log vol}_\e^{p_0}(O_p)=
-  \delta_X \hbox { log vol}_\e^\prime
 (O_p ) = -\half\int_\e^{+\infty} \hbox{tr}^\prime
[\delta_X(\tau_p^*\tau_p)
 e^{-t \tau_{p_0}^*\tau_{p_0}}]dt\eqno (3.7)$$
\s If the maps $X\mapsto  \delta_X \hbox { log vol}_\e
 (O_p )$  and $X\mapsto \delta_X \hbox{log vol}_\e^{p_0}(O_p)$ are bounded
linear maps on
the closure $H_p$ of $T_p\PC$ for the norm  induced by $<\cdot, \cdot>_p$,
 then
the preregularised mean
 curvature vector $\SC^\e$ is a vector in $H_p$ defined by $ \hbox{tr} H_X^\e=
<\SC^\e,X>_p$.\m
\item{2)} If the Riemannian structure on $\GC$ is independent of $p$,
the orbit of any point $p_0$   is a preregularisable orbit and
$$\hbox{tr} \H_{  X }^\e =
-  \delta_X \hbox { log vol}_\e^\prime
 (O_p ) = -\half\int_\e^{+\infty} \hbox{tr}^\prime
[\delta_X(\tau_p^*\tau_p)
 e^{-t \tau_{p_0}^*\tau_{p_0}}]dt\eqno (3.7 bis)$$
  where  $\hbox{tr}^\prime$ means we have restricted to the orthogonal of
  the kernel of $\tau_{p_0}^*\tau_{p_0}$ and $vol_\e^\prime$ means that
 only consider eigenvalues $\l_n^p$ that
 correspond to eigenvectors that do not belong to
 $\IC_{p_0}^p \hbox{Ker} \tau_{p_0}^*\tau_{p_o}$.
   \m {\bf Remarks}:
 \item{1)} In finite dimensions, for a compact connected
 Lie group acting via isometries on a Riemannian manifold $\PC$ of dimension
$d$,   we have for any $\e>0$ and using the
 various definitions of the volumes:
$$\eqalign{\lim_{\e \to 0}\delta_X\hbox{log vol}_\e (O_p)
 &=
  \delta_X\hbox{  log vol}_{reg} O_p   \cr
&=\delta_X\hbox{  log Vol}_{reg}  O_p\cr
&=\delta_X \hbox{ log vol}  O_p\cr} $$
 Hence going to the limit $\e \to 0$ on either side of  (3.7 bis)    we find:
 $$\hbox{tr} \HC_X = -\delta_X\hbox{ log vol} O_p.$$
If the G\^ateaux-differentiability involved is a
 $C^1$- G\^ateaux-differentiability, this yields
$$\hbox{tr} S^p= -\hbox{grad log vol} O_p$$
 This leads to    a well known result, namely
   (Hsiang's theorem [H]) that
 the orbits of $G$ whose volume are extremal among nearby orbits
 is a minimal submanifold of  $M$.
\item{2)} Equality (3.7) tells us
 that  whenever the Riemannian structure on $\G$ is independent of $p$ (as in
the case of
 Yang-Mills theory),
strongly minimal orbits of a preregularisable principal fibre bundle
are    pre-extremal
  orbits.     This  gives a weak (in the sense that we only get a sufficient
condition for strong minimality and not for minimality) infinite
dimensional version
of Hsiang's [H] theorem.
\item{3)} If both the spectrum of $\tau_p^*\tau_p$   and the Riemannian
structure on
$\G$  are
 independent of
$p$, as in the case of Yang-Mills theory in the abelian case (where the
spectrum only depends on a fixed Riemannian structure
on the manifold $M$), the orbits are strongly minimal (see also [MRT] par.5).
   \m
{\bf Proof of Proposition 3.2}:   We set $B_p=\tau_p^*\tau_p$. for the sake of
simplicity, we assume that $ B_p$ is injective on its domain, the general case
then easily follows.
\item{1)} From
 the preregularisability of the principal bundle  follows (see proposition 2.1)
 that the map $p\mapsto \hbox{det}_\e(B_p)$ is G\^ateaux-differentiable
 in the direction $X$ at point $p_0$ and
 $$  \delta_X \hbox{log det }_\e(B_p)= \int_\e^{+\infty}dt
 \hbox{tr} (\delta_X B_p e^{-tB_p} )
$$ On the other hand, by lemma 3.1
$$   \half< \int_\e^{+\infty}dt
(\delta_X B_p e^{-tB_p} )T_n^p, T_n^p>_p -e^{-\e \l_n^{p_0}}(\delta_X \IC^p
T_n^{p_0},
T_n^{p_0})_{p_0}= -
<\HC_{X }^\e
 \tilde U_n^p,
   \tilde U_n^p>_p\eqno(*)$$
The fibre bundle being preregularisable, by the results of proposition 2.0,
 the  first term on the left
 hand side is the general
 term of an absolutely convergent series. On the other hand,
the orbit being of type $(\TC)$,  the series with general term
$e^{-\e\l_n^{p_0}} (\delta_X \IC^p T_n^{p_0}, T_n^{p_0})_{p_0}$ is also
absolutely
convergent.
Hence the right hand side of (*) is absolutely convergent and $\HC_X^\e$ is
trace class since
 $(\tilde U_n)_{n\in \N}$ is a complete orthonormal basis of Im$\tau_p$.
$$-   \int_\e^{+\infty}dt
\hbox{tr}(\delta_X B_p e^{-\e B_p} )  = \hbox{tr}\HC^\e_{X }-
\delta_X \hbox{log Vol}_\e^{p_0}(B_p)
=-
 \delta_X \hbox{log det }_\e(B_p)$$
which then yields (3.7).\s
 The second part of point 1) of the proposition  follows  from the
definition of the mean
  curvature vector
in the case of a Hilbert manifold.
\item{2)} This follows from  the above  and point 3) of lemma 3.1 and holds
for any orbit $O_p$ of a regularisable fibre bundle since it does not
involve $\delta_X \IC_p$.
  \kasten
  \b  The following proposition gives an interpretation of
$\hbox{tr}_{reg}H_X$ in
  terms of the variation of
  $\hbox{vol}_{reg}
(O_p )$.
\b {\bf Proposition 3.3}: The fibres of a regularisable
 principal
 fibre bundle with structure group equipped with a fixed
 (p-independent) Riemannian metric  are heat-kernel and zeta-function
 regularisable. \s
 \item{1)} Orbits  are   heat-kernel  minimal whenever  they are
heat-kernel  extremal.
 \s
More precisely, for any point
$p_0 \in \PC$ and any horizontal
vector $X$ at point $p_0$,
 The one parameter family $\HC_X^\e$ has  a limit trace
$\hbox{tr}_{reg} \HC_X$ and
       $$\eqalign{&\hbox{tr}_{reg} \HC_{ X }= -  \delta_X
\hbox{ log vol}_{reg}( O_p) \cr
&=  \half \left[ \sum_{j=-J,j\neq 0}^{m-1} \delta_X{  b_j(p) \over j}+\int_0^1
{\delta_X  F_p(t)\over t}dt+\int_1^\infty dt  t^{-1}
 \delta_X\hbox{ tr}e^{-t\tau_{p }^* \tau_{p }} \right]
 \cr}\eqno(3.8)$$
If  the coefficients $b_j(p)$ are extremal at point $p_0$ and if $\HC_X$ is
 trace-class, then
$\hbox{tr}\HC_X=-\delta_X \hbox{log vol}_{reg} O_p$.
\item{} $$\hbox{Tr}_{reg} \HC_X=- \delta_X \hbox{ log vol}_{reg}
( O_p). \eqno(3.9)$$
  \item{2)}
The  operator $\HC_X$ has a well defined zeta function regularised trace
and  the orbit is zeta function  minimal if and only if it is
 zeta function extremal. More precisely, we have:
 $$\eqalign{\hbox{Tr}_{reg} \HC_X &=  -  \delta_X \hbox{   log Vol}_{ reg}
(O_p) \cr
 & = \hbox{tr}_{reg} \HC_X +\half \gamma   \delta_X b_0 \cr}
\eqno (3.10)$$
   \item{} If moreover $\delta_X b_0 =0$ for any
 horizontal vector $X$ at point $p_0$,  an orbit is  heat-kernel   minimal
    whenever it is zeta function minimal.
\item{3)}Whenever the map $X\mapsto  \hbox{tr}_{reg} H_X $
is  a bounded linear map on $H_p$ (with the notations of proposition 2.2),
the   regularised mean curvature vector ${\bf S}_{reg}$ is
 a vector in $H_p$
defined by $<{\bf S}_{reg},X>_p= \hbox{tr}_{reg} H_X $.
   \m {\bf Remarks:}\item{1)} In the case of a
compact connected
 Lie group acting via isometries on a finite dimensional
Riemannian manifold $\PC$ of dimension $d$,the two notions of minimality
coincide  since $b_0=d$, $\hbox{Vol}_{reg}(O_p)=\hbox{vol} (O_p)$
   and (1.10) yields:
 $$\hbox{tr}S^p=  -  \hbox{grad}\hbox{   log vol}
(O_p)  $$
where $S^p$ is the second fundamental form.
    It tells us   that
  the orbits of $G$, the volume of which  are extremal among nearby orbits is a
minimal submanifold of  $\PC$. This  proposition therefore gives an   infinite
dimensional version
of Hsiang's [H] theorem.
\item{2)}
   A zeta function formulation of Hsiang's theorem in infinite dimensions
 was already discussed in      [MRT1] in the context of
   Yang-Mill's theory. However,  there was an obstruction due to the factor
$b_0(p)$
in the zeta-function regularisation procedure which does not appear here since
 it has been taken care of in definition  (3.5) (see also [MRT2]). A
formula similar to
(3.10) (but using zeta function regularisation) can be found in [GP]
(see
formula (3.17) combined with formula (A.3)).
 \item{3)}    Proposition 3.3 puts zeta function  regularisability  and
heat-kernel regularisability on the same footing,
showing that for regularisable principal fibre bundles defined by an
isometric group action
both notions of
 regularisability hold.
\item{}It also shows that the two notions of minimalitydo not coincide in
general,
since they differ by a local term $grad b_0$, they coincide
 whenever $b_0$ is independent of $p$.
 \m
{\bf Proof of Proposition 3.3}: As before, we set $B_p=\tau_p^*\tau_p$.
As before, we shall assume for simplicity that $B_p$ is injective; the
proof then
easily extends to the case when the dimension of the kernel is locally constant
on each geodesic containg $p_0$.
\item{1)}
Since the fibre bundle is regularisable, we know by Proposition 2.1 that
the map $p\mapsto \hbox{det}_{reg}(
B_p)$ is G\^ateaux-differentiable in the direction $X$.
Let us now check that $\HC_X^\e$ has a regularized limit trace,
applying lemma 1.0.
 For this, we first investigate the
  differentiability of the map
$\e \mapsto \hbox{tr}\HC_X^\e$. By the result of Proposition 3.2,  we have
$$\hbox{tr} \HC_X^\e=  \half \int_\e^\infty dt \delta_X
\hbox{ tr}{ e^{-tB_p}
 \over t}  =-\half \delta_X \hbox{  log det}_\e(B_p) $$
The differentiability in $\e$ easily  follows from the shape of
the middle expression.
\s
 Setting as before
$F_p(t)=\hbox{tr} e^{-t B_p} -\sum_{j=-J }^{ m-1}   b_jt^{j\over m}$,
 we have
 furthermore
  $$ \eqalign{ {\partial \over \partial \e} \hbox{tr} \HC_{ X }^\e &=
-\half \e^{-1} \delta_X\hbox{  tr} e^{-\e B_p} \cr
 &= -\half \delta_X{  F_p (\e)\over \e} -
\half \sum_{j=-J }^{m -1} \delta_X b_j \e^{j-m\over m}.\cr}
$$
{}From the regularisability of the fibre bundle follows that
  $\vert\delta_X   { F_p(\e)\over \e} \vert \leq K$ for
 some $K>0$ and
$0<\e<1$
 (see assumption (2.5 bis)) which in turn implies that
 $${\partial \over \partial \e} \hbox{tr} \HC_{ X }^\e
  \simeq_0 -\half  \sum_{j=-J-m}^{-1} \delta_X b_{j+m}
 \e^{j\over m}.
  $$
 Setting $A_\e\equiv \HC_{  X  }^\e$ in Lemma 1.0, we can   define the
 regularised limit trace (replacing $J$ by $J+m$ and $a_j$ by $-\half \delta_X
 b_{j+m}$)
$$\eqalign{\hbox{tr}_{reg}\HC_{  X }
&=\lim_{\e \to 0} (\hbox{tr} \HC_{  X }^\e +
 \half \sum_{j=-J}^{-1} m { \delta_Xb_j\over j} \e^{j\over m} +
 \half \delta_Xb_0 \hbox{log} \e)\cr
&= \lim_{\e \to 0} -\half \left( \delta_X\hbox{ log det}_\e
 (B_p) - \sum_{j=-J}^{-1}
m\delta_X{ b_j\over j} \e^{j \over m}  -  \delta_X b_0  \hbox{log}
 \e\right) \quad \hbox{ by (3.7 bis )}\cr
 &= \lim_{\e \to 0} -\half\delta_X\left( \hbox{log det}_\e
  (B_p)-
\sum_{j=-J}^{-1} {mb_j \over j}
\e^{j\over m} -b_0 \hbox{log} \e \right) \quad   \cr
 &= -\half \delta_X\hbox{  log det }_{reg} (B_p) \quad \hbox{ by (1.6)}
\cr
&= \half [ \sum_{j=-J,j\neq 0}^{m-1} {m\delta_X b_j\over j}+ \int_0^1 dt
{\delta_X F_p(t) \over t}+\int_1^\infty  t^{-1}
\delta_X \hbox{tr} e^{-t\tau_p^*\tau_p}dt] \quad \hbox{by (2.6.b)}\cr}$$
 \item{} When the coefficients $b_j$ are extremal at $p_0$, we have
$$\hbox{tr} \HC_X^\e= \hbox{tr} \HC_X^\e+\half
\sum_{j=-J}^{-1}   {m\delta_Xb_j \over j}
 \e^{j\over m}+ \half \delta_X b_0 \hbox{log} \e $$
so that $\lim_{\e \to 0}\hbox{tr} \HC_X^\e= \hbox{tr}_{reg}\HC_X.$
 Since $\HC_X$ is trace class by assumption,
 going back to the definition of $\HC_X^\e$, one sees that
$\lim_{\e \to 0}\hbox{tr} \HC_X^\e=\hbox{tr} \HC_X$, which yields the
second point in 1)
of proposition 3.3.
  \item{2)}
It is well known that  the expression $  \Gamma(s)^{-1} \int_0^\infty t^{s-1}
\sum_n e^{-t\l_n}$ is finite for  Re$s$ large enough and that it has a
meromorphic continuation to the whole plane.
 Since $\Gamma(s)= (s-1) \Gamma(s-1)$, we have for $s$ with large
enough real part:
$$ \eqalign{ &\Gamma(s)^{-1} \int_0^\infty t^{s-1} \sum_n e^{-t\l_n^p}
\delta_X \l_n^p dt  =
  (s-1)^{-1}{1\over \Gamma(s-1)} \int_0^\infty t^{s-1} \sum_n e^{-t\l_n^p}
\delta_X \l_n^p dt\quad   \cr
 &=-(s-1)^{-1}{1\over \Gamma(s-1)}\int_0^\infty
 t^{s-2 }
\delta_X\hbox{ tr} e^{-tB_p}dt\quad   \hbox{see  assumption (2.2) and lemma
3.1} \cr
&=-(s-1)^{-1}{1\over \Gamma(s-1)}\left( \sum_{j=-J}^{m-1} \int_0^1
 t^{{j\over m}+s-2} \delta_X b_j  dt\right.+ \cr
&\left.+\int_1^\infty t^{s-2}
 \delta_X\hbox{  tr} e^{-tB_p}  dt +
\int_0^1 \delta_X F_p(t)  t^{s-2} dt\right) \hbox{ by (2.5)}\cr
&=-(s-1)^{-1}{1\over \Gamma(s-1)}\left[ \sum_{j=-J}^{m-1}
{1\over {j\over m}+s-1} \delta_Xb_j \right.\cr
&\left. +
\int_1^\infty t^{s-2} \delta_X\hbox{   tr } e^{-tB_p} dt +
 \int_0^1 t^{s-2} \delta_X F_p(t)  dt\right]
\cr}
$$
setting
 $F_p(t)=\hbox{tr}e^{-\e B_p}-\sum_{j=-J}^{m-1}    b_j(p)
t^{j\over m} $.
Hence, since $\Gamma(s)^{-1}= s+\gamma s^2+O(s^3)$  around $s=0$,
  going to the limit $s\to 1$, we find:
$$\eqalign{&\lim_{s\to 1}[
 \Gamma(s)^{-1} \int_0^\infty t^{s-1} \sum_n e^{-t\l_n^p}
 \delta_X \l_n^p dt  + (s-1)^{-1} \delta_X b_0 (p)] =\cr
 &=\lim_{s\to 0} (-1-\gamma s+O(s^2  )) \left[ \sum_{j=-J, j\neq 0}^{m-1}
{1\over {j\over m}+s} \delta_X  b_j  +
\int_1^\infty t^{s-1} \delta_X\hbox{  tr } e^{-tB_p}  dt\right. \cr
& \left. +
 \int_0^1 t^{s-1}\delta_X F_p(t)  dt-\gamma \delta_Xb_0\right] \cr
&=\delta_X \hbox{  det}_{reg} (B_p) -\gamma \delta_Xb_0
\quad \hbox{ by formula (1.6)  and (2.6 b)}\cr
&= -2\hbox{tr }_{reg} \HC_X-\gamma \delta_Xb_0\cr
&= \delta_X\hbox{log Det}_{reg}(B)\quad \hbox{ by formula (1.10 bis)}\cr}$$
where    $\lim_{s\to 0}\int_1^\infty t^{s-1} \delta_X\hbox{  tr } e^{-tB_p}
dt=
\int_1^\infty t^{ -1} \delta_X\hbox{  tr } e^{-tB_p}  dt$  holds
   using estimate (*) arising in the proof of proposition 2.0 and
$\lim_{s\to 0}\int_0^1 t^{s-1}\delta_X F_p(t)  dt +s^{-1} \delta_Xb_0
=\int_0^1 t^{ -1}\delta_X F_p(t)  dt
 $ by (2.4 b) and using dominated convergence.\s
 The rest of the assertions of 2) then easily follow.
\item{3)} This is a consequence of 2) using   the definition of the
regularised mean principal curvature vector
 in a Hilbert manifold.\kasten\b
    \vfill \eject \noindent
{\bf Appendix 0: G\^ateaux-differentiability}\b
We  extend here the classical notion of G\^ateaux-differentiability
 on Hilbert spaces to  Hilbert manifolds. We refer the reader to
      [AMR]  for example for the case of Banach spaces and    [E] for a
 strong version
of this notion on infinite dimensional differentiable manifolds.\s
Let $\PC$ denote  a Hilbert manifold modelled on a Hilbert space $K$ and
equipped with a (possibly   weak) Riemannian structure, which
 induces on $T_p\PC$
 a scalar product denoted by $<\cdot, \cdot>_p$. Denote by $H_p$ the
closure of $T_p\PC$ for the norm induced by $<\cdot, \cdot>_p$.
 We shall assume that
 this Riemannian structure induces an exponential map and a connection with
the usual
 properties. In particular
this yields a local diffeomorphism from the tangent bundle to the manifold.
\m
{\bf G\^ateaux-differentiability of a function }\s
 A function $
 p\to f(p)\in \R$ is said to be {\it G\^ateaux-differentiable} in the
direction $ X $ where
$ X$ is a  vector at point
$p_0\in \PC$ if $$ \delta_Xf \equiv
\lim_{\kappa \to 0}\kappa^{-1}(f( p_\kappa)-f(p_0) )$$
 exists   with
 $p_\kappa= \hbox{exp}_{p_0} \kappa   X  $, exp denoting the
geodesic exponential on $\PC$.
A function $
 p\to f(p)\in \R$ is said to be   G\^ateaux-differentiable  at a point $p_0$
if it is
 G\^ateaux-differentiable
in all directions at that point.
 \s The map $X\to \delta_X f$ defined on $H_p$ is denoted by
 $\delta_{p_0}f$. If $f$ is  G\^ateaux
 differentiable   at point
 $p_0$
 and if
 the map $ \delta_{p_0}f$  is a bounded linear form   on
  $H_{p_0}$,
 by Riesz theorem, we can  identify it with
  a vector $G_{p_0}f$  in $H_{p_0}$:
 $$<G_{p_0} f,X>_{p_0}= \delta_Xf \quad \f X\in T_{p_0} \PC$$ If the
 Riemannian structure is strong, $G_{p_0}f$ lies in
$T_{p_0} \PC$.\s
If $f$ is differentiable at point $p_0$, $f$ is G\^ateaux-differentiable
 at $p_0$ and
$  \delta_{p_0}f= d_{p_0}f$ is a bounded linear map on $T_p\PC$.
 Hence, if the Riemannian structure is strong,
$ d_{p_0}f$ is identified to a vector $\hbox{grad}_{p_0}f\in T_{p_0}\PC$ by
 $d_{p_0}f(X)=<\hbox{grad}_{p_0}f,X>_{p_0}$ and we have
 $\hbox{grad}_{p_0}f = G_{p_0}f$.
  \b If for any $p_0\in \PC$ and for any vector $X$ at $p_0$,
the map $  \delta_{p_0}f $ is a bounded linear form on
 the Hilbert space $H_p$ and if for any point $p_0$, the map
 $p \mapsto   \delta_pf\circ \hbox{exp}_{p_0}^*$ is continuous on
exponential charts at $p_0$
for the operator norm, the function $f$ will be called
  $C^1$-G\^ateaux-differentiable. Although G\^ateaux-differentiability is
weaker than differentiability,   $C^1$ G\^ateaux-differentiability implies
 $C^1$ differentiability.
\m
  {\bf Lemma  0.1}: A   $C^1$-G\^ateaux-differentiable function is $C^1$
differentiable.
\m {\bf Proof}: The proof goes as in the case of a Hilbert space (see
 [AMR] Corollary 2.4.10) using the exponential map. \s
  Let $V(p_0)$ be a neighborhood of $p_0$ small enough so as to be in the
image of the exponential map
at point $p_0$ and on which $p\mapsto G_pf$ is uniformly continuous.
 For   $p\in V(p_0)$, there is a vector
$X\in T_{p_0}\PC$ such that
$\hbox{exp}_{p_0}X=p$ and we shall set
 $p_\kappa= \hbox{exp}_{p_0}\kappa X$.  Then
$$\eqalign{\vert f(p )-f(p_0)- \delta_{p_0}fX\vert&=\vert
\int_0^1( {d\over d\kappa}f(p_\kappa)-  \delta_{p_0}X) d\kappa\vert\cr
&\leq \int_0^1  \vert   \delta_{p_\kappa}f\hbox{exp}_{p_0*}(  \kappa X) (
X)-   \delta_{p_0}fX\vert d\kappa\cr
&\leq \int_0^1  \Vert\vert
\delta_{p_\kappa}f\hbox{exp}_{p_0*}(\kappa \cdot)(\cdot)  -  \delta_{p_0}f\Vert
 \vert
\Vert X\Vert d\kappa\cr
&\leq \hbox{sup}_{\kappa\in [0,1]}   \Vert\vert
\delta_{p_\kappa}f\hbox{exp}_{p_0*}(\kappa \cdot)(\cdot)  -
 \delta_{p_0}f\Vert
 \vert
\Vert X\Vert  \cr}$$
 and the r.h.s converges to zero as $X$ goes to $0$.
  This says that  $df (p_0)$ exists and is equal
to $  G_{p_0}f$ so that $f$ is differentiable and $C^1$.\kasten
 \m
{\bf Remark}: In finite dimensions, one can show in a similar way that
$C^1$-G\^ateaux-differentiability implies $C^1$-differentiability.
\b
{\bf G\^ateaux-differentiability of operators}\s
 For a family of bounded  operators
$B_p $ on a given Hilbert space $H$, $p$ varying in $ \PC$,
     we shall say that the map $p\to   B_p$
is
   G\^ateaux   differentiable  in  the direction $X$ at
 point $p_0$
if
 $ \kappa^{-1}( B_{p_\kappa}-B_{p_0})  $
converges  as $\kappa$ goes to zero in norm to an operator $\delta_XB$.
 \m A family of unbounded operators $B_p$ defined on a common dense domain
$D$ in a given Hilbert space
 $H$ to a given Hilbert space $K$, will be called G\^ateaux-
differentiable at point $p_0$ in the direction $X$ if
for every $u$ in  $D$,
$ \kappa^{-1}   B_{p_\kappa}u-B_{p_0}u $ converges to a vector in $K$ when
$\kappa \to 0$, thus defining a densily defined operator $\delta_XB$
on $H$.
Furthermore, if $p\mapsto B_p^*$ is G\^ateaux-differentiable and if
the spaces $H$ and $K$ are equipped with fixed scalar products, then
  $\delta_X B_p^*\subset(\delta_X B_p)^*$ . Indeed,   for any
$u $ in the domain
of $\delta_X B_p$, $v$ in the domain of $\delta_X B_p^*$, we have:
$$\eqalign{ <\delta_X B_p u,v>_K &= \lim_{\kappa \to 0}
< {B_{p_\kappa}-B_{p_0} \over \kappa}u,v>_K\cr
&=  \lim_{\kappa \to 0} <u, {B_{p_\kappa}^*-B_{p_0}^* \over \kappa}v>\cr
&=<u, \delta_X B_p^*v>\cr}$$
In particular,
the G\^ateaux-differential at point $p_0$
of  a family of rotations $ R_p$ of $H$ which coincide with identity at
 point $p_0$
is antisymmetric. Indeed, differentiating the relation
$R_p^*R_p=I$ at point $p_0$, we obtain
$\delta_XR_p^*=-\delta_X R_p$. Since $\delta_X R_p^*\subset
(\delta_X R_p)^*$,
we find that $-\delta_XR_p\subset (\delta_X R_p)^*$, which says that $R_p$ is
antisymmetric.
\b {\bf Covariant G\^ateaux-differentiability }
  \s Let $O(\PC)$  be the orthonormal bundle   and let us associate to
any $C^1$ curve $\sigma$ in $\PC$
 its horizontal lift $\tilde \sigma$ in $O(\PC)$ given $\tilde \sigma(0)$.
Thus, to    a curve
$\sigma(\kappa)$, we associate
$\tilde \sigma(\kappa)\in \hbox{Isom}(H, T_{p_\kappa}\PC)$ where $p_\kappa=
\sigma(\kappa)$.
\s For $p_0\in \PC$ and $X\in T_{p_0} \PC$, define the curve
$  \sigma(\kappa)= \hbox{exp}_{p_0} \kappa X$.\m
 {\bf Definition}: Let $p_0 \in \PC$ and $X$ be a vector at point $p_0$.
A vector field $V$ on $\PC$ will be called
G\^ateaux-differentiable  at point $p_0$ in the direction $X$ if, setting
$\sigma(\kappa)= \hbox{exp}_{p_0}\kappa X=p_\kappa$,
  if the expression
  $
  \kappa^{-1}\tilde \sigma(0)
 (\tilde \sigma(\kappa)^{-1} V(p_\kappa)-\tilde \sigma(0)^{-1} V(p_0))
 $
has a  well defined limit in $T_{p_0} \PC$ when $\kappa$ goes to zero. This
limit is
denoted by $\delta_X V$ so that
 $$\delta_XV= \tilde \sigma(0)
 {d\over d\kappa}_{\kappa=0}(\tilde \sigma(\kappa)^{-1} V(p_\kappa)).
 $$
\m Let $p\mapsto \TC_p$ be a  field of operators  acting from a dense
 space
$D\subset H$ of a fixed Hilbert space $H$ into $T_p\PC$.
The map $\TC_p$ will be called G\^ateaux-differentiable at point $p_0$
 in the direction $X$ if for every
vector $u\in D$, the vector field $\TC_pu$ is G\^ateaux- differentiable
 at point $p_0$ in the direction $X$. We shall write
 $\delta_X  ( \TC_pu)=\delta_X \TC_p u$.  \m
 \vfill \eject
     { \bf  Appendix A. Regularizability  and minimality  of the orbits for
a loop group coadjoint
action }\m
In this paragraph, we want to investigate the heat-kernel regularizability
and minimality of the orbits
 for the coadjoint action of the loop group in agreement with
 the work by King and Terng [KT], who showed they are
 zeta function regularisable and minimal submanifolds.   We shall take the
same notations as in [KT]. \s
 Let $G$ be a connected, compact Lie group, ${\bf g}$ its Lie algebra,
$(\cdot, \cdot)_0$ a fixed
 $Ad$ invariant inner product on ${\bf g}$. Let
$$\O \equiv \{ g\in H^1([0, 1], G), g(0)=g(1)=e\}$$
With the notations of the previous paragraph, we take $\G=\O$, which
is a Hilbert Lie group. The exponential map is given by
$(\hbox{Exp} \xi)(t)= \hbox{exp} \xi(t)$.We take
$\PC=\HO$.\s
$\O$ acts on $\HO$ on the left by
$$ \eqalign{ \O\times \HO &\to \HO\cr
(g, \gamma)&\to g\gamma g^{-1}-g^\prime g^{-1}= g\star \gamma.\cr}$$
 The action is free, smooth and isometric (see [KT] and references therein).
 The orbit space is $G$ and the map
$$\eqalign {\pi:\HO&\to G\cr
\gamma &\to  g (1)\quad \hbox{with} \quad  g^{-1}  g^\prime=
\gamma, \quad   g (0)=e\cr}$$
yields a principal fibre bundle structure of $\HO$ with structure
 group $\O$, the fibres of which are congruent w.r.to isometries of
 $\HO$. For $g\in H^1([0,1], {\bf g})$, $\gamma\in H^1([0,1],   G )$ with
 $g(0)=e$, we have
 $\pi(\gamma  \star g)= \pi(\gamma)g(1)^{-1}$. Hence for another $\tilde
 \gamma\in H^1([0,1], {\bf g})$, choosing $g(1)= \pi(\tilde \gamma)^{-1}
\pi(\gamma)$,
 we have
 $\pi(\gamma\star g)=\pi(\tilde \gamma)$.
 \s
In the following, $\hat a$ denotes the constant loop with fixed value $a$.\s
 For fixed $\gamma\in H^0([0,1], {\bf g})$, the tangent map at unit element
$\hat e\in \O$ to the map:
$$\eqalign{ \theta_\gamma:\O&\to \HO\cr
g&\to g\star\gamma\cr}$$
  is given by
$$\eqalign{\tau_\gamma:\OO&\to\HO \cr
u &\to [ u, \gamma]-u^\prime\cr}\eqno (A.1)$$
 where $\OO\equiv T_{\hat e} (\O)= \{ u\in H^1([0,1], {\bf g}),
 u(0)=u(1)=0\}
$.\s Notice that $\tau_\gamma$ is the sum of a first order differential
operator
 independent of $\gamma$ and of a bounded operator which depends on
$\gamma$. Of course,
  if $G$ is abelian the $\gamma$-dependent part
 vanishes.
\m
 Let us equip $\GC$ with a fixed weak
right invariant Riemanian metric defined  by the inner product
$$<u,v>=\int_0^1 (u(t), v(t))_0 dt$$
on the space $$T_g(\O)=\{u\cdot g, u\in H^1([0, 1],g), u(0)=u(1)=0\}$$ by
$$<u_1\cdot g,  u_2\cdot g>= \int_0^1 (u_1(t), u_2(t))_0dt$$
where $u\cdot g= R_{g*} u$.\s
 The adjoint $\tau_\gamma^*$ of $\tau_\gamma$ w.r.to this scalar product is
given by  $\tau_\gamma^*(\xi)=[  \gamma,\xi ]+ \xi^\prime$ for $\xi  \in
H^1([0, 1],g), \xi(0)=\xi(1)
$  and $\tau_\gamma^*\tau_\gamma$ is self adjoint on
$$D(\tau_\gamma^*\tau_\gamma)\equiv \{  u\in H^2([0, 1],g), u(0)=u(1),
 \quad u^\prime(0)=u^\prime(1)\}.$$  Notice
that it is not injective on its domain since it contains the constant loops.
\s  \b {\bf Proposition A.1}: The orbits are heat-kernel preregularizable
and
strongly minimal in $\HO$.
 The preregularized volumes of the orbits are constant.
\m {\bf Proof}:
 \s Since every orbit contains a constant loop (see [Se] Proposition 8.2 )
and since the action is isometric, it is enough to consider the orbits of
 constant loops. Since
the orbits are congruent via isometries
(see [KT]), is also enough show the heat-kernel regularisability and strong
minimality of
 the orbit containing the 0-loop, namely the orbit of  $\gamma_0=\hat 0$.\s
 Since Ker$\tau_{\hat 0}^*=\{\hat a, a\in {\bf g}\}$,  the normal space to
the zero orbit at point $\hat 0$ is the space of
constant loops.
We   therefore want  to show that for any $\e>0$ the operaotr $\HC_{\hat x}^\e$
is trace-class for $\e>0$ and $\hbox{tr}\HC_{\hat x}^\e=0$
for any constant loop $\hat a$, $\hat x$.
\s In order to apply proposition 2.2, we first check the
 preregularisability of the bundle. Frorm the above considerations,
follows that it is enough to check
assumptions (2.1)-(2.3) for $p_0=\hat 0$.
 Let $x\in g$  and let  ${\bf t}$ be the maximal abelian subalgebra of
${\bf g}$
 containing $x$. The constant loop $\hat x$ will play the role of the
tangent  vector
 $X$ at point $p_0$ used in the main bulk of the paper. Let $a=sx$, $s\in
[0,1]$, so that $\hat a$ lies on the geodesic starting at $\hat 0$ in
direction $\hat x$.
The complexified Lie algebra ${\bf{ g_{\C}}}$ has an orthonormal basis
built up from
$z_{ \alpha}, \alpha \in \Delta^+$,
and $a_k, k=1, \cdots, r$, $r$ being the dimension of ${\bf t}_{\C}$  and
$\Delta^+$
  the set of positive roots of $G$.
 They satisfy
 the following anticommuting  relations:
$$[a,z_\alpha]= -i \alpha(a) z_\alpha.$$
Setting
$z_\alpha=x_\alpha+i y_\alpha$, this yields:
$$[a,x_\alpha]=\alpha(a) y_{\alpha}, [a,y_\alpha]=
-\alpha(a) x_\alpha, [a,a_i]=0.\eqno (A.2)$$
\s Let
$  {\O}_{\C} $ denote the  complexification  of $\O$  and $\overline
{{\O}_{\C}}$ its closure
w.r. to the hermitian form $H(\cdot, \cdot)$ induced by the
$L^2$ scalar product $<\cdot, \cdot>  $.
 An orthonormal basis of $\overline{ {\O}_{\C}}$ is given by
$$z_{n, \alpha}(t)= z_\alpha e^{2i\pi n t}, a_{k, n}(t)
 =
 a_k e^{2i\pi n t}, n\in \Z, \alpha \in \Delta^+, k=1,\cdots, r.$$
This yields an orthonormal   basis of the  (real) $L^2$ closure
 $\overline{ \O}$ of $\O$:
 $$\eqalign{ r_{\alpha, n}&= \sqrt 2\hbox{Re} (z_{  \alpha,n}) ,
s_{\alpha, n}= \sqrt 2
\hbox{ Im}(z_{ \alpha,n}) ,\cr
  \rho_{k, n}&= \sqrt 2\hbox{Re} (a_{k,n}) ,
 \sigma_{k,n}=\sqrt 2 \hbox{Im} (a_{k,n}) ( n\neq 0) \cr}$$
\m
An easy computation using definition (A.1) and relations (A.2) shows that
 for
$u=\sum_{\alpha, n}
 u_{\alpha, n} z_{\alpha,n}+ \sum_{k, n} \eta_{k,n} a_{k,n} \in
 {{\O}_{\C}}$ $$\tau_{\hat a} (u)=
 \sum_{\alpha, n} ( i\alpha(a) -2i\pi n)u_{\alpha,n}z_{\alpha, n}-
\sum_{k,n} (2i\pi n) \eta_{k,n}a_{k,n}$$
which yields:
$$\tau_{\hat a} (r_{\alpha, n})=(2\pi n-\alpha(a))s_{\alpha, n},\quad
\tau_{\hat a}^*\tau_{\hat a}(r_{\alpha, n})=(2\pi n-\alpha(a))^2
(r_{\alpha, n})$$
$$\tau_{\hat a} (s_{\alpha, n})=(-2\pi n+\alpha(a))r_{\alpha, n},\quad
\tau_{\hat a}^*\tau_{\hat a}(s_{\alpha, n})=(2\pi n-\alpha(a))^2
(s_{\alpha, n})$$
$$\tau_{\hat a} (\rho_{k, n})=-2\pi n \sigma_{k, n},\quad
\tau_{\hat a}^*\tau_{\hat a}(\rho_{k, n})=(2\pi n )^2
(\rho_{k, n})$$
 $$\tau_{\hat a} (\sigma_{k, n})= 2\pi n \rho_{k, n},
\quad \tau_{\hat a}^*\tau_{\hat a}(\sigma_{k, n})=(2\pi n )^2 (\sigma_{k, n}),
(n\neq 0).$$
   This proves that the  operator $\tau_{ \hat a}^* \tau_{ \hat a}$
  has purely discrete spectrum
 given by $$\eqalign{ &  (2\pi n-\alpha(a))^2,
( 2\pi n +\alpha(a))^2, n\in \N^* (\hbox{ each  taken    with multiplicity 2
}),\cr
&(\alpha(a))^2 , \alpha \in \Delta^+ (\hbox{   taken    with multiplicity 2
}),  \cr
 & (2\pi n)^2,  n\in \N^*  (\hbox{  each taken    with multiplicity 4r}),\cr
&  0
 (\hbox{  taken    with multiplicity  r})
  \cr}$$
 \m  Clearly the eigemvalues beahve asymptotically as $ (2\pi  n)^2$    so
that the first assumption for a preregularisable bundle is satisfied.
Each of these eigenvalues is clearly differentiable as a function of $a$
and we  set
$$  \beta_n^{\e , \hat x}(t)\equiv
\delta_{\hat x}  (\l_n^{\hat a}) e^{-t \l_n^{\hat 0}}$$
For
$\l_n^{\hat a} = (2\pi n-\alpha(a))^2$ or  $\l_n^{\hat a} =
 (2\pi n+\alpha(a))^2$, we have
$\beta_n^{\e, \hat x}= 4\pi n \alpha(x) e^{-t(2\pi n)^2}$. For
$\l_n^{\hat a} = (2\pi n )^2$,  $\beta_n^{\e, \hat x}=0$. Hence
$\beta_n^{\e, \hat x}$ which is independent of $a$ satisfies
$\sum_n \vert \beta_n^{\e, \hat x}\vert <\infty$ which proves (2.1)
 and (2.2)
 then follows easily. The bundle is therefore preregularisable.
  \s From lemma 3.1 and the fact that $\G$ is equipped with a fixed
Riemannian structure
 then follows  that the spectrum of $\HC_{\hat x}^\e$ is
 purely discrete and given by
$$\eqalign{(\mu_n^\e, n\in \N)\equiv &\left({\alpha(x)\over 2\pi n }
e^{-\e (2\pi n)^2}, {-\alpha(x)\over 2\pi n} e^{-\e (2\pi n)^2},
 \alpha \in \Delta^+,n\in \N^*   \right.\cr
&\left. 0, \cdots, 0 (\hbox{ r  times}), n\in \N^*\right)\cr}$$
 which coincides with the expression obtained by [KT].\s
  Applying proposition 3.2, we have that  $\HC_{ \hat x}^\e$ is trace
 class, $a\mapsto \hbox{det}_\e(\tau_{\hat a}^* \tau_{\hat a})$ is
G\^ateaux-differentiable in the direction $\hat x$ at $\hat 0$ and
 $$\hbox{tr} \HC_{\hat x}^\e=
-\half \delta_{\hat x}
\hbox{log vol}_\e^\prime(O_{\hat a}) \eqno (*)$$
where $O_{\hat a}$ is the orbit of $\hat a$ and $vol_\e^\prime$ is to be
understood in the sense
of proposition 3.2.\s
Since the l.h.s of $(*)$ is the general term of an absolutely convergent
series, we can reorder the terms and write
$$\delta_X  h_\e(\tau_{\hat a}^*\tau_{\hat a}) =
2\sum_{n\in \N}(\mu_n^{\e, \hat x}+\mu_{-n}^{\e, \hat x})=0$$
so that $$\hbox{tr}\HC_{\hat x}^\e= <\gld_\e  (\tau_{\hat a}^*
\tau_{\hat a}), \hat x>_{\hat 0}=0$$
  Thus the orbits are heat minimal  the preregularized volumes of the
orbits are constant.\kasten
{\bf Remark}: This confirms the fact proved by King and Terng [KT]
that the orbits are zeta function regularizable and minimal and that the
zeta function regularized volume of the
fibres are constant. \b
 {\bf Appendix B. Gauge orbits for Yang-Mills theory}\b
In this appendix, we confront the notions of
 regularity and minimality
 based on
heat-kernel regularizations
 developped above with those developped in [MRT1]
 based on zeta function regularizations
   in the case of Yang-Mills theory. This
     appendix does not offer new results, its
     only purpose being to relate
     the general framework developped in the bulk
     of this paper with the concrete example of Yang-Mills theory.
 This appendix is  loosely written and we refer
 the reader to [MRT1] for a  precise presentation
 of the case of Yang-Mills orbits.\s

 The comparison of these two approaches shows that even
in this particular example, they are not equivalent in general. Only in
very particular cases does the notion of heat-kernel
minimality coincide with that of zeta function minimality.
\m
For the description and results  that follow, we  refer the reader  to
 [KR], [FU]
 and [MV]  (and many other  references mentioned therein)
where  the gauge action in Yang-Mills theory was investigated in details.
 \s Let $\pi:P\to P/G\simeq M$ be a smooth
 principal bundle over an $n$-dimensional smooth oriented manifold $M$,
 where $G$ and $M$ are compact connected, $G$ being a Lie
 group with Lie algebra ${\bf g}$. We shall specialise $G$ to
a matrix group here.
\s  If we allow manifolds with boundary,
we recover the framework  dual to the one described in section II
 (in the sense that the action on the left on  ${\bf g}$ valued
vector fields there is replaced here by an action on the right on
 ${\bf g}$- valued forms), when applying
 the setting we are about to describe to $M=[0,1]$ with Dirichlet boundary
conditions.
However, for the sake of clarity, we take $M$ boundaryless here and have
treated the
dimension 1 case with boundary separately in section II.\m
 Automorphisms of $P$ are described
 equivalently
\item{1)} as fibre preserving homeomorphisms $\phi:P\to P$
 commmuting with the $G$-action
\item{2)} as continuous sections  of the principal fibre bundle
 $E_G\equiv P\times_G G$ associated to $P$, given by the quotient
of $P\times G$
by the right action
$$ \eqalign{(P\times G)\times G &\to P\times G\cr
(p,g_1)\cdot g &\to (p\cdot g, g^{-1} g_1 g)\cr}$$
which is free and proper.
\item{3)}as continuous functions $\phi:P\to G$ satisfying
$$\phi(p\cdot g)=g^{-1} \phi(p) g$$
A connection  on $P$ can be described equivalently  as
\s  1)   a subbundle  $A\subset TP$ such that
\item{a)} $\pi_*A_p=T_{\pi(p)}M$ and $A_p\cap T_pF_{\pi(p)}=\{0\}$,
\item{b)}$A_{p\cdot g}=(R_g)_*A_p$, for $g\in G, p\in P$.
\s where $F_x$ is the fibre over $x$.\m
 2) a ${\bf g}$-valued one form $A$ on $P$ such that
 \item{a)} $\f X\in {\bf g}, \f p \in P$, $A(\tau_p X)=X$
\item{b)}$ \f g\in G, \f v\in T_pP$, $A_{p\cdot g}(R_g)_*v=
 Ad g^{-1} A_p(v)$.
 \m The action of automorphisms on connections can be described as follows.
An automorphism $\phi$ of class $C^2$ acts on a connection of class
$C^1$ in the following way:
$$(\phi\cdot A)_{\phi(p)}=T_p \phi\cdot A_p\eqno (B.1)$$
Expressing $A$ in terms of ${\bf g}$-valued one forms and $\phi$
 as a $G$ valued function on $P$, we have
$$(\phi\cdot A)_p=Ad \phi(p) \circ A_p-R_{\phi(p)*}^{-1}\circ
 T_p\phi$$
It is clear from this formula that $\phi\cdot A=A\quad \f  A$ if and
 only if
$\phi$ is constant and $Ad \phi=Id$, i.e $\phi\in Z(G)$, the center of
$G$.
\b Let us, with the notations of section I,  describe the group ${\bf G}$.
 Let $k>{n\over 2}+1$ be an integer and let ${\bf G^k}$ denote the
automorphisms of $P$ of class $H^k$, i.e the  group of
 $H^k$ sections of the bundle $E_G$, which is a Hilbert Lie group
modelled on $H^k(E_{\bf g})$,
 the space of $H^k$ sections of the adjoint bundle $ E_{\bf g}$ which is
the vector bundle associated to
$P$ with fibre ${\bf g}$ defined as the quotient of $P\times {\bf g}$
for the action of the group $G$  given by:
 $$\eqalign{(P\times {\bf g})  \times G &\to E_{\bf g}\cr
              ((p, \xi),g)&\to (p\cdot g, ad g \xi)\cr}$$
$H^k(E_{\bf g})$ can also be seen as the space of $H^k$
  $\bf g$ valued
 equivariant functions on $P$:
$$H^k(E_{\bf g})=\{\phi\in H^k(P,{\bf g}), \f g\in G, \f p\in P,
 \phi(p\cdot g)= g^{-1} \phi(p) g\}$$
 Let now $\overline{\bf G^k}$ denote the quotient of ${\bf G^k}$ by its center.
It is also a Hilbert lie group modelled on $H^k(E_{\bf g})$. With the
notations of section I, we set
$\bf G\equiv \overline{\bf G^{k+1}}$.
\b Let us now, with the notations of section I,  describe the manifold
$\PC$ on which ${\bf G}$ acts.
The space $\AC^k$ of connections of class $H^k$ is a closed affine subspace of
the Hilbert space
$H^k( {\bf g}\otimes T^*P)$ of   $H^k$.  By the Sobolev embedding
theorems,
$H^k$ connections are of class $C^1$ and $H^{k+1}$ automorphisms of $P$ of
class $C^2$.
Hence $\overline{\bf G^{k+1}}$ acts on $\AC^k$ and we
  shall call a connection $A\in \AC^k$
 irreducible when $\{ \phi\in \overline{\bf G^k}, \phi\cdot A= A\} $ is
reduced to the
 constant map with value $1$. Let $\overline{\AC^k}$ denote the space of
irreducible connections in $\AC^k$. It is a Hilbert manifold modelled on
$H^k( {\bf g}\otimes T^*P)$. We shall set $\PC=\overline{\AC^k}$.
\b We shall consider the following action of $\G$ on $\PC$:
$$\eqalign{ \overline {\G^{k+1}} \times \bar \AC^k&\to \AC^k\cr
(\phi, A)&\to \phi\cdot A\cr}\eqno (B.2)$$
where $\phi\cdot A$ is given by (B.1).
 \b One can show that for $k>{n \over 2}+3$, the action of ${\bf G} $
on $\PC$ is free, smooth
 and proper and that the canonical projection
$\PC\to \PC/{\G}$ yields a principal fibre bundle structure on the moduli
space $\PC/ \G $.
\m The space $\overline{\AC^k}$ can be equipped with a $\G$-invariant $L^2$
(hence weak)
Riemannian structure given by
$$<\alpha, \beta>\equiv\int_P (\alpha(p), \beta(p) )d\mu(p)$$ where
 $\mu$ is a smooth $G$-invariant measure on $P$ (recall that $G$ is compact),
 $\alpha, \beta$ ${\bf g}$
valued one forms on
$P$ and $(\cdot, \cdot)$ the bundle metric in ${\bf g} \otimes T^*P$
induced by an AdG invariant scalar product on
${\bf g}$ and a $G$ invariant Riemannian metric on $P$.
The group $\G$ is equipped with a fixed  $L^2$ $\G$-invariant
structure in a similar way. The Riemannian structure on $\overline \AC^k$
induces a Riemannian connection (see e.g [KR]).
 \b  For given $A\in {\overline{A^k}}$, the tangent map at the constant map
with value $1_G$ to
$$\eqalign{\theta_A : {\overline {\bf G^{k+1}}}&\to {\bar \AC^k}\cr
     \phi&\to \phi\cdot A\cr}$$ reads
$$\eqalign{\tau_{A }: H^{k+1} (E_{\bf g})&\to H^k({\bf g}\otimes T^*P)
\cr
\l &\to ad\l \circ A-T \l=-\nabla_A\l\cr}$$
where $\nabla_A$ is the covariant derivative in the associated bundle
$E_{\bf g}$ defined by the connection $A $ on $P$ or equivalently
$\nabla_A=\nabla+[A,\cdot]$
where $\nabla$ is the covariant derivative ${\bf g}\to {\bf g} \otimes
T^*\PC$ defined by the Riemannian metric on $\PC$.
Notice that taking $M=S^1$, we recover formula (A.1 ).
 Notice that when $G$ is abelian, $\tau_A =-\nabla$ so that $\tau_A$ is
independent of the choice of the connection $A$.
It is easy to check that for an  irrreducible
connection $A$, the map $\tau_A$ is injective.\s
   For a $\Ci$ connection, if  $\tau_A^*$ denotes the adjoint
of $\tau_A$ w.r. to the $L^2$ structure  on $\PC$ and on $\G$,
 the operator $\tau_A^* \tau_A$ is an elliptic operator of order 2 with
 $\Ci$ coefficients (see e.g [KR]).
It is by now a well known result that
$\tau_A^*\tau_A$ satisfies  assumptions
 (2.1)-(2.5) (see e.g [MRT1])
 so that an orbit $O_{A_0}$ for the action (B.2) is
   regularisable given  the assumption that
 the dimension of the kernel of $\tau_A^*\tau_A$
   is locally constant around $A_0$. This
   last assumption is discussed in [MRT1]. Hence the following conclusions:
    \m
\item{\bf 1)}The orbit of any smooth irreducible connection  is regularisable.
 It is  minimal if and only if it has minimal   heat kernel
regularised volume   among orbits of the same type.
 \item{ } If the group $G$ is abelian, the orbit of any smooth
 irreducible connection  is  (strongly) minimal.\s
  The  proof of this statement in the
  general case  follows from proposition 3.3.
The result in the abelian case  follows from proposition 3.2 using the fact
that $\tau_A^*\tau_A$
  is independent of $A$ and hence so are its eigenvalues
 since they only involve the Riemannian structures on $P$ and $g$ but not
connections.
\m
  \item{\bf 2)} If the dimension of $M$ is odd or equals 2,
  the two notions of minimality (heat-kernel an zeta function minimality)
coincide and if
  the dimension of $M$ is 4, this holds for  smooth
 irreducible Yang-Mills connections.
   \s  To prove this statement,
we check that  the  G\^ateaux differential in
any horizontal direction $X$
of the coefficient $b_0$
 in the   heat-kernel expansion
 vanishes which then yields the result by proposition 3.3.
\item{} By [G], we know that it vanishes when the dimension of $M$ is odd and
 that when dim$M$=2,
   $b_0(A)= c_1\int_M s(g) + (\hbox{dim }Ad P)(\hbox{Vol} M)$
  where
$c_1$ is a constant, $s(g)$ the scalar curvature of a metric $g$ on $M$
so that it is independent of $A$ and $\delta_Xb_0=0$.\item{}
 When dim$M$=4, $b_0(A)= c_2(g)+c_3\int_M \hbox{tr} \vert
 F_A\vert^2$, where $F_A$ is the curvature of $A$,
where $c_2(g) $ only depends on the metric $g$ and $c_3$ is a nonzero constant.
Since Yang-Mills connections are exactly the ones   which extremise the
 Yang-Mills functional $\int_M \vert F_A\vert^2$, we have
$\delta_X b_0=
\delta_X\int_M \vert F_A\vert^2=0$. (see [MRT1]).
 \m
       \vfill \break \noindent
 {\bf REFERENCES}
\item{[AMT]} R.Abraham,J.E. Marsden, T.Ratiu, {\it Manifolds, tensor
analysis and applications}
{\it Global Analysis, Pure and Applied}, Modern Methods for the study of
non linear phenomena in engineering,
Addison Wesley, 1983
\item{[AP1]} M.Arnaudon, S.Paycha,
{\it Factorization of semi-martingales on infinite dimensional
 principal bundles}, to appear in Stochastics and Stochastic Reports
\item{[AP2]} M.Arnaudon, S.Paycha, {\it The geometric and physical
relevence of some stochastic tools on Hilbert manifolds}, manuscript
 \item{[AJPS]} S.Albeverio, J.Jost, S.Paycha, S.Scarlatti,
{\it A mathematical introduction to string \hfill \break \noindent
theory-variational problems, geometric and probabilistic methods}, to
appear  ( Cambridge University Press)
 \item{[FU]} D.S Freed, K.K. Uhlenbeck, {\it Instantons and Four manifolds},
Springer Verlag (1984)
\item{[FT]} A.E.Fischer, A.J. Tromba, {\it On a purely Riemannian proof of
the structure and dimension of the ramified moduli space of a compact
Riemann surface}, Mathematische Annalen {\bf 267}, p.311-345, (1984)
\item{[G]} P.B.Gilkey, {\it Invariance Theory, The heat equation and the
Atiyah-Singer index theorem}, Publish or Perish, Wilmington, 1984
\item{[GP]}:D.Groisser, T.Parker, {\it Semi-classical Yang-Mills theory
I:Instantons}, Comm.Math.Phys. {\bf  135}, p.101-140 (1990)
\item{[H]}W.Y. Hsiang, {\it On compact homogeneous minimal submanifolds},
Proc.Nat. Acad. Sci. USA {\bf 56} (1966) p5-6
  \item{[KR]} W.Kondracki, J.Rogulski, {\it On the stratification of the
orbit space},
Dissertatione Mathematicae Polish Acad. of Sci.{\bf 250} (1986), p. 1-62
  \item{[KT] }C.King, C.L.Terng, { \it Volume and minimality of submanifolds
 in path space}, in  "Global Analysis and Modern Mathematics", ed. K.
Uhlenbeck,
Publish or Perish ( 1994)
\item{[MRT1] }Y.Maeda, S.Rosenberg,P. Tondeur,
{\it The mean curvature of gauge orbits},  in  "Global Analysis and Modern
Mathematics", ed. K. Uhlenbeck,
Publish or Perish ( 1994)
\item{[MRT2] }Y.Maeda, S.Rosenberg,P. Tondeur,
{\it Minimal orbits of metrics and Elliptic  Operators}, manuscript 95
 \item{[MV]} P.K.Mitter, C.M.Viallet, {\it On the bundle of connections and
the gauge
 orbit manifold in Yang-Mills theory}, Comm.Math. Phys. 79 (1981), P. 457-472
 \item{[P]} S.Paycha {\it Elliptic operators in the functional quantisation for
gauge field theories},   Comm.Math.Phys  {\bf 166}, p.433-455 ( 1995)
\item{[RS]} D.S.Ray, I.M.Singer {\it R-torsion and the Laplacian on
Riemannian manifodls},
Advances in Mathematics 7, p.145-210 (1971)
\item{[S]} G.Segal {\it Unitary representations of some infinite
Dimensional groups}, Comm. Math. Phys. {\bf 80}, p. 301-342 (1981)
 \end